\documentclass[11pt]{article}
\usepackage[utf8]{inputenc}
\usepackage[T1]{fontenc}
\usepackage{comment}
\usepackage{authblk}
\usepackage{relsize}
\usepackage{changepage}
\usepackage{mathtools}
\usepackage[english]{babel}
\usepackage{amsmath}
\usepackage{amsthm}
\usepackage{amsfonts}
\usepackage{amssymb}
\usepackage{bbm}
\usepackage{graphicx}
\usepackage[usenames,dvipsnames]{color}
\usepackage{soul,xcolor}

\theoremstyle{plain}
\newtheorem{theorem}{Theorem}[section]

\newtheorem{lemma}[theorem]{Lemma}
\newtheorem{claim}[theorem]{Claim}

\newtheorem{algorithm}[theorem]{Algorithm}
\newtheorem{conjecture}[theorem]{Conjecture}

\makeatletter
\newcommand{\vast}{\bBigg@{4}}
\newcommand{\Vast}{\bBigg@{5}}
\makeatother

\theoremstyle{definition}

\newtheorem{definition}[theorem]{Definition}

\newtheorem{def/prop}[theorem]{Definition/Proposition}

\usepackage[left=2cm,right=2cm,top=2cm,bottom=2cm]{geometry}
\usepackage{amssymb}
\usepackage{hyperref}
\makeatletter
\def\namedlabel#1#2{\begingroup
    #2%
    \def\@currentlabel{#2}%
    \phantomsection\label{#1}\endgroup
}
\usepackage[normalem]{ulem}
\usepackage{dsfont}
\usepackage{accents}
\usepackage{verbatim}
\usepackage{ulem}
\usepackage{pgf,tikz,pgfplots}
\pgfplotsset{compat = 1.16}
\usepackage[all]{xy} 

\usepackage{stackengine}
\stackMath
\newcommand\tsup[2][2]{%
 \def\useanchorwidth{T}%
  \ifnum#1>1%
    \stackon[-.5pt]{\tsup[\numexpr#1-1\relax]{#2}}{\scriptscriptstyle\sim}%
  \else%
    \stackon[.5pt]{#2}{\scriptscriptstyle\sim}%
  \fi%
}

\DeclareMathOperator{\Prob}{\mathbb{P}}

\DeclareMathOperator{\eps}{\varepsilon}
\DeclareMathOperator{\cA}{\mathcal{A}}
\DeclareMathOperator{\cB}{\mathcal{B}}
\DeclareMathOperator{\cC}{\mathcal{C}}
\DeclareMathOperator{\cD}{\mathcal{D}}
\DeclareMathOperator{\cE}{\mathcal{E}}
\DeclareMathOperator{\cF}{\mathcal{F}}
\DeclareMathOperator{\cH}{\mathcal{H}}
\DeclareMathOperator{\cT}{\mathcal{T}}
\DeclareMathOperator{\cS}{\mathcal{S}}
\DeclareMathOperator{\cP}{\mathcal{P}}

\usepackage{subfig}

\usepackage{enumerate}

\title{\scshape
  The semi-random tree process}

\author[1]{Sofiya Burova}
\author[2]{Lyuben Lichev}
\affil[1]{Ecole Normale Supérieure de Lyon, Lyon, France}
\affil[2]{Univesit\'e Jean Monnet and Institut Camille Jordan, Saint-Etienne, France}

\begin{document}
\setstcolor{red}
\maketitle
 
\begin{abstract}
The online semi-random graph process is a one-player game which starts with the empty graph on $n$ vertices. 
At every round, a player (called Builder) is presented with a vertex $v$ chosen uniformly at random and independently from previous rounds, and constructs an edge of their choice that is incident to $v$. 
Inspired by recent advances on the semi-random graph process, we define a family of generalised online semi-random models. 

We analyse a particular instance that shares similar features with the original semi-random graph process and determine the hitting times of the classical graph properties minimum degree $k$, $k$-connectivity, containment of a perfect matching, a Hamiltonian cycle and an $H$-factor for a fixed graph $H$ possessing an additional tree-like property. Along the way, we derive a few consequences of the famous Aldous-Broder algorithm that may be of independent interest.
\end{abstract}

\section{Introduction}

The online semi-random graph process, introduced by Michaeli and studied recently in~\cite{BMPR21, B-EGHK20, B-EHKPSS20, GKMP20, GH20}, is a stochastic one-player game. It starts with the empty graph on $n$ vertices and, at every round, a player (called Builder) is presented with a vertex $v$ chosen uniformly at random and independently from previous rounds. Then, Builder constructs one of the edges incident to $v$ of their choice and proceeds to the next round.

From the grounding paper on the existence of thresholds by Bollob\'as and Thomason~\cite{BT87} to the recent breakthrough on the Kahn-Kalai conjecture by Park and Pham~\cite{PP22}, thresholds have been a major topic in the theory of random graphs. The notion of threshold has its natural analogue in settings where the random process is controlled by an agent. More particularly, many natural graph properties and their hitting times in the semi-random graph process by an ``optimal'' strategy were studied. For example, constructing a graph of minimum degree $k\in \mathbb N$ was studied in~\cite{B-EHKPSS20} where the authors proved that Builder needs $(h_k+o(1))n$ rounds to complete the task whp where $(h_k)_{k\ge 1}$ is a recursively defined sequence. Moreover, whp Builder needs between $0.9326 n$ and $1.2053 n$ rounds to construct a perfect matching~\cite{GMP21}, and analogous bounds are known for the construction of a Hamiltonian cycle: in fact, in~\cite{B-EHKPSS20} the authors showed that Builder needs between $(\log 2 + \log(1+\log 2))n$ and $3n$ rounds to construct a Hamiltonian cycle whp, and these bounds were improved in~\cite{GKMP20} to $(\log 2 + \log(1+\log 2)+10^{-8})n$ and $2.6114 n$, respectively. In fact, for every $\eps > 0$, any graph on at most $n$ vertices and sufficiently large maximum degree $\Delta$ may be constructed in $\left(\tfrac{3\Delta}{2}+\eps\right)n$ rounds whp~\cite{B-EGHK20}.

Inspired by the above model, we define a family of online semi-random hypergraph models as follows. 
\begin{definition}[The general online semi-random hypergraph process]\label{def SR general}
Fix a set $S$ of $N$ initially unmarked vertices, a family $\cF$ of subsets of $S$ and a probability distribution $\mu$ on $\cF$. At round $i\in \mathbb N$, Builder is given a set of $\cF$ sampled according to $\mu$ independently from previous rounds and they can mark any vertex in $\cF$ (including one already marked at a previous round). 
\end{definition}

Note that the original semi-random model is a particular case of this more general version, where $N = \tbinom{n}{2}$, $S$ is the set of edges in the complete graph $K_n$, $\cF = \{\{ij: j\in [n]\setminus \{i\}\}: i\in [n]\}$ and $\mu$ is the uniform probability measure on $\cF$. In general, the aim of Builder is to mark all vertices in some family of sets $\cH\subseteq 2^S$. In the above setup, this could be the set of all perfect matchings of $K_n$, the set of all Hamiltonian cycles of $K_n$, etc. We believe that the generality of Definition~\ref{def SR general} leaves little hope for proving very precise results without any restrictions on $\cF$, $\mu$ or $\cH$. Therefore, in this paper, we choose to conduct a thorough analysis of a particular instance, which we call \emph{the uniform spanning tree semi-random model} (or \emph{the USTSR model} for short) that bears close resemblance to the original semi-random graph model. Again, our set $S$ will consist of the $\tbinom{n}{2}$ vertex pairs of $K_n$, but this time $\cF$ will be the set of all spanning trees of $K_n$, while $\mu$ remains the uniform probability measure on $\cF$. This choice is motivated by the fact that, on the one hand, in both models, Builder may construct a spanning tree of $K_n$ in $n-1$ rounds with probability 1, and on the other hand, the structure of the edges in $\cF$ is somehow preserved: in both cases, these are spanning trees of $K_n$. However, note that the number of elements in $\cF$ in our case is a lot larger: there are exactly $n^{n-2}$ spanning trees of $K_n$ by Cayley's formula while only $n$ spanning stars. We remark that the construction of a fixed graph in the USTSR model was studied in the master thesis of the first author~\cite{Sof21}.

\paragraph{Notations.} We use the standard asymptotic notations $o, \omega, O, \Omega, \Theta$. Our default limit variable is $n$; when this is not the case, we indicate the limit variables in an index to the asymptotic notation, for example $O_x$ or $\Theta_{x,y}$. Moreover, for $\Theta_{x,y}$ and similar, the limit over the latter variable (in this case $y$) is taken first. We write $f = O_{\Prob}(g)$ (resp.\ $f = \Omega_{\Prob}(g)$) when $|f|\le Cg$ for some positive random variable $C$ which is a.s.\ finite (resp.\ when $|f|\ge cg$ for some positive random variable $c$ which is a.s.\ non-zero); as usual, $f=\Theta_{\Prob}(g)$ means that $f = O_{\Prob}(g)$ and $f = \Omega_{\Prob}(g)$. Also, for a sequence of probability spaces $(\Omega_n, \cF_n, \Prob_n)_{n\ge 1}$ and a sequence of events $(A_n)_{n\ge 1}$ satisfying $A_n\in \cF_n$ for all $n\ge 1$, we say that $(A_n)_{n\ge 1}$ happens \emph{with high probability}, or \emph{whp}, if $\Prob_n(A_n)\to 1$ as $n\to \infty$.

For a graph $G = (V(G),E(G))$, we call $|G| := |V(G)|$ the \emph{order} of $G$ and $|E(G)|$ the \emph{size} of $G$. For a vertex $v\in V(G)$ and a set $U\subseteq V(G)$, we denote $N_U(v)$ the neighbourhood of $v$ in $U$ (i.e. the set of all neighbours of $v$ in $U$), and $N_U[v] = N_U(v)\cup \{v\}$. We further denote $V_n = V(K_n) = \{v_i\}_{i=1}^n$ and $E_n = E(K_n)$. The graphs $(K_n)_{n\ge 1}$ will be considered embedded in each other so in particular the sequence $(V_n)_{n\ge 1}$ is increasing with respect to inclusion. For a given graph $H$, an $H$-factor on $n$ vertices is a vertex-disjoint collection of $\lfloor |H|^{-1}n\rfloor$ copies of $H$.\footnote{Note that most sources require $|H|$ to divide $n$ in the definition of an $H$-factor; this additional condition will not be of much importance in our setting.}

For a set $S$ and any integer $k\ge 0$, we denote by $\tbinom{S}{k}$ the set of subsets of $S$ of size $k$ and also $2^S = \bigcup_{k\ge 0} \tbinom{S}{k}$. We also set $\cT_n$ to be the family of spanning trees of $K_n$.

\paragraph{Formal description of the process.} A \emph{deterministic strategy} $\sigma$ is a mapping from $2^{E_n}\times \mathcal T_n$ to $E_n$ such that for all $F\in 2^{E_n}$ and $T\in \mathcal T_n$ one has $\sigma (F,T)\in E(T)$ (here, $F$ corresponds to the set of edges already added at previous stages of the USTSR process, $T$ corresponds to the spanning tree of $K_n$, sampled at the present round, and $\sigma(F,T)$ is the edge of $T$ that Builder chooses to construct). Note that it is possible that $\sigma(F,T)$ is already in $F$, particularly when $E(T)\subseteq F$. A \emph{randomized strategy} is a random variable on the set of deterministic strategies. We denote by $\mathcal S_n$ the set of randomized strategies for our USTSR process on $n$ vertices. Here and below, for all $i\ge 0$, $T_i$ is the uniform spanning tree sampled at round $i$ of the USTSR process, and $G_{\sigma, i}$ (or simply $G_i$) is the (random) graph, constructed by Builder after round $i$.

We remark that, by abuse of terminology, first, we often describe the strategy $\sigma$ that Builder adopts by providing an associated algorithm rather than a mapping, and second, we often identify ourselves with Builder when constructing and analyzing $\sigma$ and refer to Builder as ``we''.

\subsection{Our results}
A \emph{property} $\cP$ is a subset of $2^{E_n}$. For a strategy $\sigma\in \cS_n$ and a property $\cP$, we denote by $\tau(\sigma, \cP) = \tau_n(\sigma,\cP)$ the random variable $\min\{i\ge 0: G_{\sigma, i}\in \cP\}$.

Let $\cP_{\min, k}$ be the family of graphs on $n$ vertices with minimum degree $k$. The following result shows that Builder may construct a graph with minimum degree $k$ in slightly more than $kn/2$ rounds and also quantifies the error term.

\begin{theorem}\label{thm min degree}
\hfill
\begin{enumerate}[(i)]
    \item\label{pt 1 thm deg} For any $k = k(n)$ satisfying $k\le n-1$ and any strategy $\sigma \in \mathcal S_n$, whp 
\begin{equation*}
\tau(\sigma, \mathcal P_{\min, k})\ge \tfrac{kn}{2} + \tfrac{(2\pi n)^{1/2}}{4} + o(n^{1/2}).
\end{equation*}
\item\label{pt 2 thm deg} If $k = o(n^{1/2})$, then there is a strategy $\sigma\in \mathcal S_n$ such that whp
\begin{equation*}
\tau(\sigma, \mathcal P_{\min, k})\le \tfrac{kn}{2} + \tfrac{(2\pi n)^{1/2}}{4} + o(n^{1/2}),
\end{equation*}
and moreover $G_{\tau(\sigma, \mathcal P_{\min, k})}$ contains only vertices of degrees $k$ and $k+1$.
\item\label{pt 3 thm deg} Furthermore, if $k = o(n)$, then there is a strategy $\sigma\in \mathcal S_n$ such that whp
\begin{equation*}
\tau(\sigma, \mathcal P_{\min, k})\le \tfrac{kn}{2} + O(\max\{k, n^{1/2}\}),
\end{equation*}
and moreover $G_{\tau(\sigma, \mathcal P_{\min, k})}$ contains only vertices of degrees $k$ and $k+1$.
\end{enumerate}
\end{theorem}

We believe that the conclusion of Theorem~\ref{thm min degree}~\eqref{pt 2 thm deg} holds more generally for $k = o(n)$. Since, first, we do not prove this in the current work, and second, the proofs of Theorem~\ref{thm min degree}~\eqref{pt 2 thm deg} and~\eqref{pt 3 thm deg} follow the same lines, we only sketch the proof of Theorem~\ref{thm min degree}~\eqref{pt 3 thm deg} based on the ideas from the proofs of~\eqref{pt 1 thm deg} and~\eqref{pt 2 thm deg}.

\begin{conjecture}
The conclusion of Theorem~\ref{thm min degree}~\eqref{pt 2 thm deg} holds for $k = o(n)$.
\end{conjecture}

The next theorem deals with the construction of factors. For a given graph $H$, we denote by $\cP_H$ the family of graphs on $n$ vertices that contain an $H$-factor. Clearly, for all $\sigma\in \mathcal S_n$ one has $\tau(\sigma, \cP_H)\ge \frac{n |E(H)|}{|H|} - O(1)$ since every $H$-factor on $n$ vertices contains that number of edges. Also, we say that a graph $H$ is \emph{central} if there is a vertex $v\in V(H)$ (which we call a \emph{witness}) such that, for all edges $e\in E(H)$ incident to $v$, one has that $H\setminus e$ is a disconnected graph.
\begin{theorem}\label{thm factors}
\hfill
\begin{enumerate}[(i)]
    \item\label{pt 1f} For any fixed tree $H$ and any function $\omega = \omega(n)\to \infty$ there is a strategy $\sigma\in \mathcal S_n$ satisfying $\tau(\sigma, \cP_H)\le \frac{n(|H|-1)}{|H|} + \omega n^{3/4}$ whp.
    \item\label{pt 2f} For any fixed central graph $H$ there is a constant $C_H > 0$ and a strategy $\sigma\in \mathcal S_n$ such that $\tau(\sigma, \cP_H)\le C_H n$ whp.
    \item\label{pt 3f} Fix $r\ge 3$. There is a central graph $H$ containing $K_r$ for which there is a strategy $\sigma\in \mathcal S_n$ such that $\tau(\sigma, \cP_{H})\le \frac{n|E(H)|}{|H|} + \omega n^{3/4}$ whp.
\end{enumerate}
\end{theorem}

Note that the statement of Theorem~\ref{thm factors}~\eqref{pt 3f} may seem a bit surprising at first sight since even factors of graphs, containing very dense subgraphs ($K_r$ in particular), may have construction time which, up to lower order terms, is the same as the easy deterministic lower bound given just before Theorem~\ref{thm factors}. In our example $K_r$ is only a ``small part'' of the graph $H$ so only ``a few'' copies of $K_r$ have to be constructed to see an $H$-factor appear. Of course, many ``almost complete'' copies of $K_r$ are also present at that stage. In order to avoid ``wasting'' the edges used for their construction, we make sure that $H$ contains ``many'' copies of almost complete graphs on $r$ vertices as well.

An obvious weakness of Theorem~\ref{thm factors} is that it does not treat certain natural graphs as, for example, cycles and complete graphs. In fact, for any fixed graph $H$, fixing $\lfloor n/|H|\rfloor$ predetermined copies of $H$ and waiting until all of them are constructed, takes $O(n\log n)$ rounds whp. Despite the fact that this naive strategy could possibly be far from optimal, we believe that the true optimal construction time is ``somewhere in the middle''.
\begin{conjecture}
Let $H$ be a fixed non-central graph. Then, for every strategy $\sigma\in \cS_n$, $\tau(\sigma, \cP_H)\gg n$ whp.
\end{conjecture}

Denote by $\cP'_H$ the set of graphs, in which every vertex participates in a copy of $H$. Note that the property $\cP'_H$ is a relaxation of $\cP_H$ since, up to restricting $n$ to multiples of $|H|$, $\cP_H\subseteq \cP'_H$. Also, denote by $d_{\min}(H)$ the minimum degree of $H$.

\begin{theorem}\label{thm covering}
For every graph $H$ there is a strategy $\sigma\in \cS_n$ such that $\tau(\sigma, \cP'_H)\le d_{\min}(H) n + o(n)$ whp. 
\end{theorem}

Note that the upper bound in Theorem~\ref{thm covering} is optimal up to a factor of 2: indeed, if every vertex is covered by a copy of $H$, then the degree of every vertex must be at least $d_{\min}(H)$. \\ 

\noindent \textbf{Remark.} More careful analysis of our approach may be conducted to ensure that the $o(n)$ error term may be replaced by $n^{1-\eps}$ for some $\eps = \eps(H) = \Omega(|E(H)|^{-1})$. We choose to spare the details for the sake of a less elaborate presentation of the proof. \\

The next theorem concerns the appearance of a Hamiltonian cycle. Denote by $\cP_{HC}$ the family of graphs on $n$ vertices that contain a Hamiltonian cycle.

\begin{theorem}\label{thm HC}
Fix any function $\omega = \omega(n)\to \infty$ as $n\to \infty$. Then, there is a strategy $\sigma\in \cS_n$ satisfying whp
\begin{equation*}
    \tau(\sigma, \cP_{HC})\le n + \omega n^{3/4}.
\end{equation*}
\end{theorem}

\noindent
Another property of interest to us is $k$-connectivity. For $k\in [n-2]$, denote by $\cP_{\mathrm{con},k}$ the set of $k$-connected graphs on $n$ vertices.

\begin{theorem}\label{thm k-conn}
Fix any $k = o(n)$ and any function $\omega = \omega(n)\to \infty$. Then, there is a strategy $\sigma\in \cS_n$ satisfying whp
\begin{equation*}
\tau(\sigma, \cP_{\mathrm{con}, k})\le \tfrac{nk}{2} + \omega k \max\{k, n^{3/4}\}.
\end{equation*}
Moreover, if $k = o\left(\tfrac{n^{1/2}}{(\log n)^2}\right)$, there is a spanning $k$-connected subgraph of $G_{\tau(\sigma, \cP_{\mathrm{con}, k})}$ with vertices of degrees $k$ and $k+1$ only.
\end{theorem}

We remark that the second term in Theorem~\ref{thm k-conn} could possibly be reduced in some cases; in fact, in our proof the second terms originate from the construction of $k$ independent perfect matchings. However, a similar idea to the proof of Theorem~\ref{thm min degree}~\eqref{pt 2 thm deg} may be used to derive a smaller order second term when $k = o\large(\large(\tfrac{n}{\log n}\large)^{1/2}\large)$. Roughly speaking, one would still need to construct $k$ almost perfect matchings consecutively but give up on each of them when there are $2(n\log n)^{1/2}$ unmatched vertices left, and continue with the next matching. Nonetheless, for every new matching, one could begin the construction by connecting the vertices that failed to be connected in previous matchings with vertices that did not fail to do so (this is where the bound on $k$ is used). Then, in the end, one would have $k$ disjoint independent sets of size $2(n\log n)^{1/2}$ containing only vertices of degree $k-1$, and one could thus try to extend all $k$ perfect matchings simultaneously. This observation may be used to obtain a second term of the form $O(n)$ for $k = o\large(\large(\tfrac{n}{\log n}\large)^{1/2}\large)$; we omit the details.

\paragraph{Plan of the paper.} In Section~\ref{sec prelims} we present several preliminary results and in Section~\ref{sec proofs} we prove the theorems stated above.

\paragraph{Acknowledgements.} We are grateful to Dieter Mitsche for several related discussions and to two anonymous referees for multiple useful comments.

\section{Preliminaries}\label{sec prelims}

\subsection{Probabilistic preliminaries}

\subsubsection{Concentration inequalities}

The next lemma recalls the well-known Chernoff's bound, see e.g. Corollary 2.3 in~\cite{JLR}.
\begin{lemma}[Chernoff's bound]\label{chernoff}
Let $X \sim \mathrm{Bin}(n,p)$ be a binomial random variable with $\mathbb E[X]=np=\mu$. Then, for every $\delta\in [0,1]$,
\begin{equation*}
\mathbb P(|X-\mu| \ge \delta \mu) \le 2\exp\bigg (-\dfrac{\delta^2 \mu}{3}\bigg).
\end{equation*}
\end{lemma}

Another inequality that will be of use to us is Azuma's martingale inequality, see e.g. Theorem 2.25 in~\cite{JLR}.

\begin{lemma}[Azuma's inequality]\label{lem azuma}
Fix a martingale $X_0 = 0, X_1, \dots, X_n$ and a sequence of constants $(c_k)_{k=1}^n$ satisfying $|X_{k-1}-X_k|\le c_k$ for all $k\in [n]$. Then, for all $t > 0$, 
\begin{equation*}
\max\{\Prob(X\ge t), \Prob(X\le -t)\}\le \exp\left(-\frac{t^2}{2\sum_{k=1}^n c_k^2}\right).
\end{equation*}
\end{lemma}

\subsubsection{The differential equation method}
The differential equation method is a powerful and, by now, a classical technique used to follow the evolution of a discrete random process. Pioneered by Kurtz~\cite{K70} around 1970, the method was first applied in a combinatorial context by Wormald~\cite{Nick1, Nick2, Nick3}. The basic idea is the following:
given a sequence of discrete random variables $(X_t)_{t \ge 0}$ with bounded increments, the trajectory $(X_t)_{t \ge 0}$ (properly rescaled) is approximated by the solution of an ordinary differential equation suggested by the expected changes $(\mathbb E[X_{t+1}-X_t])_{t\ge 0}$. Tight concentration of the process $(X_t)_{t \ge 0}$ around the idealized trajectory is provided by classical results from martingale theory.

The precise formulation of the theorem given here is taken from~\cite{Lutz}: a function $f$ is said to be $L$-Lipschitz on $D \subseteq \mathbb R^{\ell}$ if $|f(x) - f(x')| \le L \max_{1 \le k \le \ell} |x_k - x'_k|$ holds for all points $x=(x_1, \ldots, x_{\ell})$ and $x'=(x'_1, \ldots, x'_{\ell})$ in $D$, where $\max_{1 \le k \le \ell} |x_k - x'_k|$ is the $\ell^{\infty}-$distance between $x$ and $x'$.
\begin{theorem}[\cite{Lutz}, Theorem 2]\label{Thm:DEMethod}
Given $a, n \ge 1$, a bounded domain $D \subseteq \mathbb{R}^{a+1}$, functions $(F_k)_{k\in [a]}$ with $F_k: D \to \mathbb{R}$, and $\sigma$-algebras ${\mathcal{F}}_{0}\subseteq {\mathcal{F}}_{1} \subseteq \ldots$, suppose that the random variables $(Y_k(i))_{k\in [a]}$ are ${\mathcal{F}}_{i}$-measurable for $i \ge 0$. Suppose also that for all $i \ge 0$ and all $k \in [a]$, the following holds whenever $\left(\tfrac{i}{n}, \tfrac{Y_1(i)}{n}, \ldots, \tfrac{Y_a(i)}{n}\right) \in D$:
\begin{enumerate}
\item\emph{Trend hypothesis:} $\left| \mathbb{E}[Y_k(i+1) - Y_k(i) \mid  {\mathcal{F}}_{i}] - F_k\left(\tfrac{i}{n}, \tfrac{Y_1(i)}{n}, \ldots, \tfrac{Y_a(i)}{n}\right)\right| \le \delta$ for some $\delta \ge 0$, with $F_k$ being $L$-Lipschitz for $L \in \mathbb{R}$.
\item\emph{Boundedness hypothesis:} $|Y_k(i+1) - Y_k(i)| \le \beta$ for some $\beta > 0$,
\item\emph{Initial condition:} $\max_{k \in [a]}|Y_k(0) - \hat{y}_kn| \le \lambda n$ for some $\lambda > 0$, for some $(0, \hat{y}_1, \ldots, \hat{y}_a) \in D$.
\end{enumerate}
Then, there are $R=R(D,(F_k)_{k\in [a]},L) \in [1, \infty)$ and $T=T(D) \in (0, \infty)$ such that, for any $\lambda \ge \delta \min \{T, \tfrac{1}{L}\} + \tfrac{R}{n}$, with probability at least $1-2a\exp\left(-\tfrac{n\lambda^2}{8T\beta^2}\right)$ we have
\begin{equation*}
\max_{i \in [\sigma n]} \max_{k \in [a]} |Y_k(i) - y_k(\tfrac{i}{n}) n| \le 3\exp(LT) \lambda n,
\end{equation*}
where $(y_k(t))_{k\in [a]}$ is the unique solution to the system of differential equations $y_k'(t)=F_k(t,y_1(t), \ldots, y_a(t))$ satisfying $y_k(0)=\hat{y}_k$ for all $k\in [a]$, and $\sigma=\sigma(\hat{y}_1, \ldots, \hat{y}_a) \in [0,T]$ is any choice of $\sigma \ge 0$ with the property that $(t, y_1(t), \ldots, y_a(t))$ has $\ell_{\infty}$-distance at least $3\exp(LT)\lambda$ from the boundary of $D$ for all $t \in [0, \sigma)$.
\end{theorem}

\subsection{The Aldous-Broder algorithm}
Perhaps the simplest algorithm for sampling a uniform spanning tree of a given (finite) graph is the Aldous-Broder algorithm~\cite{Ald90, Bro89}. 

\begin{algorithm}[Aldous-Broder algorithm]\label{AB algo}
\emph{Input: a graph $G$. Pick an arbitrary vertex of $G$ and start a simple random walk $(W_i)_{i\ge 0}$ from it. Set $E = \emptyset$. At any step $i\ge 1$, if $|E| < |V(G)|-1$, do:
\begin{itemize}
    \item if $W_i\in \{W_j: j\in \{0,\dots,i-1\}\}$, then continue;
    \item if $W_i\notin \{W_j: j\in \{0,\dots,i-1\}\}$, then do $E\leftarrow E\cup \{W_{i-1}W_i\}$ and continue.
\end{itemize}
Output: the set $E$.}
\end{algorithm}

Note that one may replace the simple random walk on $G$ with a lazy random walk given by the transition matrix
\begin{equation*}
P(u,v) = \mathds{1}_{v\in N_G[u]} \tfrac{1}{\deg_G(u)+1},
\end{equation*}
and the above algorithm would still output a uniform spanning tree of $G$. Although such replacement does not make sense from an algorithmic point of view, this observation will prove useful to us in several of the proofs below. We suspect that (some of) the following lemmas may already be present in the literature, but in the absence of a satisfactory reference, we chose to provide complete proofs.

\begin{lemma}\label{lem independence}
Fix two disjoint sets of vertices $S_1,S_2\subseteq V_n$. Given $T\sim \mathrm{Unif}(\cT_n)$, the random graphs $H_1 = (S_1, E(T)\cap \tbinom{S_1}{2})$ and $H_2 = (S_2, E(T)\cap \tbinom{S_2}{2})$ are independent of each other.
\end{lemma}
\begin{proof}
We rely on the Aldous-Broder algorithm (Algorithm~\ref{AB algo}) for generating a uniform spanning tree of $K_n$ via a lazy random walk $(W_t)_{t\ge 0}$. The statement follows from the fact that $H_1$ and $H_2$ are measurable with respect to the projection of $(W_t)_{t\ge 0}$ onto $S_1$ and $S_2$, respectively, and $S_1\cap S_2 = \emptyset$. 
\end{proof}

\begin{lemma}\label{lem clique AB}
Given $T\sim \mathrm{Unif}(\mathcal T_n)$ and any $k\le n-1$, the probability that $E(T)\cap \tbinom{V_k}{2} = \emptyset$ is $\left(1 - \tfrac{k}{n}\right)^{k-1}$.
\end{lemma}
\begin{proof}
Once more, we rely on the Aldous-Broder algorithm (Algorithm~\ref{AB algo}) for generating a uniform spanning tree of $K_n$ via a lazy random walk. Start a random walk $(W_t)_{t\ge 0}$ from a vertex in $V_n\setminus V_k$ and, for every $i\in [k]$, let

$$\tau_i = \min\{j\in \mathbb N: \{W_t: t\in \{1,\dots,j\}\} \text{ contains }i \text{ elements from }V_k\}.$$ 
Then, the event $E(T)\cap \tbinom{V_k}{2} = \emptyset$ may be rewritten as 
$$\forall i\in [k], W_{\tau_i-1}\notin \{W_{\tau_j}: j\in \{1,\dots,i-1\}\} \text{ and } W_{\tau_i+1}\notin V_k\setminus \{W_{\tau_j}: j\in \{1,\dots,i\}\}.$$

In other words, for $E(T)\cap \tbinom{V_k}{2} = \emptyset$ to take place, the following happens: every time the random walk visits a vertex in $V_k$ for the first time, it starts at a vertex in $V_n\setminus V_k$, visits $V_k\setminus (W_{\tau_j})_{j=1}^{i-1}$ and then goes back to a vertex in $(V_n\setminus V_k)\cup (W_{\tau_j})_{j=1}^{i-1}$ or stays in place. Note that, under this event, $W_{\tau_i+1}$ is a vertex in $(V_n\setminus V_k)\cup (W_{\tau_j})_{j=1}^i$ chosen uniformly at random and the event $\{W_{\tau_{i+1}-1}\notin (W_{\tau_j})_{j=1}^{i}\}$ has probability $\tfrac{n-k}{n-k+i}$. Letting by convention $\{W_{\tau_0+1}\notin V_k\}$ to be the trivial event with probability 1, we deduce that
\begin{align*}
&\Prob(E(T)\cap \tbinom{V_k}{2} = \emptyset)\\
=\hspace{0.3em} 
&\Prob(\forall i\in [k], W_{\tau_i-1}\notin (W_{\tau_j})_{j=1}^{i-1} \text{ and } W_{\tau_i+1}\notin V_k\setminus (W_{\tau_j})_{j=1}^i)\\
=\hspace{0.3em} 
&\prod_{i=1}^k \Prob(W_{\tau_i-1}\notin (W_{\tau_j})_{j=1}^{i-1} \text{ and } W_{\tau_i+1}\notin V_k\setminus (W_{\tau_j})_{j=1}^i\mid  W_{\tau_{i-1}+1}\notin V_k\setminus (W_{\tau_j})_{j=1}^{i-1})\\
=\hspace{0.3em} 
&\prod_{i=1}^k \frac{n-k}{n-k+i-1}\cdot \frac{n-k+i}{n}\\
=\hspace{0.3em} &\frac{(n-k)^k}{\tfrac{(n-1)!}{(n-k-1)!}}\cdot \frac{\tfrac{n!}{(n-k)!}}{n^k} = \frac{(n-k)^{k-1}}{n^{k-1}},
\end{align*}
which concludes the proof of the lemma.
\end{proof}

\begin{lemma}\label{lem markov chain}
Fix a Markov chain $(M_t)_{t\ge 0}$ with three states, $A,B$ and $C$, and transition probabilities 
$$P(A,B) = p, P(B,A) = q, P(A,C) = 1-p, P(B,C) = 1-q \text{ and } P(C,C) = 1.$$
Then, conditionally on $M_0 = A$, the probability that the transition $A\to C$ was used is $\tfrac{1-p}{1-pq}$ if $pq < 1$ and $1$ otherwise.
\end{lemma}
\begin{proof}
The case $pq = 1$ is clear. Let $pq < 1$. Then, all trajectories that use the transition $A\to C$ may be written as $A(BA)^*C$ where the transition from $A$ to $B$ and back may repeat itself an arbitrary number of times. Thus, the total probability is given by
\begin{equation*}
\sum_{i=0}^{\infty} (1-p)p^iq^i = \frac{1-p}{1-pq}.
\end{equation*}
\end{proof}

The next lemma computes the probability that the uniform spanning tree of $K_n$ does not have common edges with a fixed complete bipartite graph.

\begin{lemma}\label{lem 2.8}
Given $T\sim \mathrm{Unif}(\mathcal T_n)$ and any $k,\ell\in \mathbb N$ satisfying $k+\ell\le n$, 
$$\Prob(E(T)\cap (V_k\times (V_n\setminus V_{n-\ell})) = \emptyset) = \frac{(n-k)^{\ell-1}(n-\ell)^{k-1}(n-k-\ell)}{n^{k+\ell-1}}.$$
\end{lemma}
\begin{proof}
We rely on the Aldous-Broder algorithm again. This time, however, our strategy is to directly compute the number of spanning trees of the graph $H_n = K_n\setminus (V_k\times (V_n\setminus V_{n-\ell}))$. To this end, we find the probability that the spanning tree of $H_n$ generated by the lazy random walk $(W_t)_{t\ge 0}$ is given by the path $1,2,\dots,n$.

Let us start a lazy random walk $(W_j)_{j\ge 0}$ on $K_n$ from the vertex $v_1$. Then, for every $i\in [k]$, the probability that $v_{i+1}$ is the $(i+1)$-st vertex, added to the tree, conditionally on the fact that $v_i$ was the $i$-th vertex added to the tree, is given by $\tfrac{1}{n-\ell} + \tfrac{1}{(n-\ell) (n-\ell-i)}$: indeed, there is a probability of $\tfrac{1}{n-\ell}$ that the random walk visits $v_{i+1}$ immediately after its first visit of $v_i$, and otherwise it stays within $V_i$ at the next step with probability $\tfrac{i}{n-\ell}$ and the first edge going out of $V_i$ towards $V_{n-\ell}\setminus V_i$ is $v_iv_{i+1}$ with probability $\tfrac{1}{i(n-\ell-i)}$. This gives the product
\begin{equation*}
\prod_{i=1}^{k} \frac{1}{n-\ell}\left(1 + \frac{1}{n-\ell-i}\right) = \prod_{i=1}^{k} \frac{(n-\ell-i+1)}{(n-\ell)(n-\ell-i)} = \frac{1}{(n-\ell)^{k-1}(n-\ell-k)}. 
\end{equation*}

Now, for every $i\in [k+1, n-\ell]$, we will compute the probability that $v_{i+1}$ is the vertex added immediately after $v_i$ conditionally on the fact that the vertex set, explored by the random walk up to now, is $V_i$ and the vertex $v_i$ was the last one explored. There are two cases: either the random walk visits the vertex $v_{i+1}$ immediately after it visits $v_i$, which happens with probability $\tfrac{1}{n}$, or it remains within $[i]$ with probability $\tfrac{i}{n}$. In the second case, we consider the Markov chain $(M_{1,j})_{j\ge 0}$ with three states, $A_1$, $B_1$ and $C_1$, and transition probabilities
\begin{align*}
p_1 =\hspace{0.3em} &P_1(A_1,B_1) = \tfrac{k}{n-i+k},\\ 
q_1 =\hspace{0.3em} &P_1(B_1,A_1) = \tfrac{i-k}{n-\ell-k},\\ &P_1(A_1,C_1) = 1-p_1,\\ 
&P_1(B_1,C_1) = 1-q_1,\\
&P_1(C_1,C_1) = 1.    
\end{align*}

Then, by setting $A_1 = V_i\setminus V_k$, $B_1 = V_k$ and $C_1 = V_n\setminus V_i$, we construct the Markov chain, obtained from the random walk by forbidding the steps between vertices in $A_1$ and in $B_1$. Note that these transitions only make the hitting time of $C_1$ longer without influencing the probability of hitting $C_1$ from a vertex in $A_1$ or from a vertex in $B_1$. Then, by Lemma~\ref{lem markov chain}, the probability of a transition $A_1\to C_1$ when $M_{1,0} = A_1$ is
$$\frac{1-p_1}{1-p_1q_1} = \frac{1-\tfrac{k}{n-i+k}}{1-\tfrac{k}{n-i+k}\tfrac{i-k}{n-\ell-k}},$$
and the probability of a transition $A_1\to C_1$ when $M_{1,0} = B_1$ is $1 - \frac{1-q_1}{1-p_1q_1} = \frac{q_1(1-p_1)}{1-p_1q_1}$. Moreover, every edge between $A_1$ and $C_1$ has equal probability to be used for the unique transition between $A_1$ and $C_1$, which is $\tfrac{1}{(i-k)(n-i)}$. Thus, the probability that $v_iv_{i+1}$ is the next edge, added to the tree, is
\begin{align*}
&\frac{1}{n} + \frac{1}{(i-k)(n-i)}\left(\frac{i-k}{n}\frac{1-p_1}{1-p_1q_1}  + \frac{k}{n}\frac{q_1(1-p_1)}{1-p_1q_1}\right)\\
=\hspace{0.3em}
&\frac{1}{n} + \frac{1}{(i-k)(n-i)} \frac{1-p_1}{n(1-p_1q_1)} \left(i-k + kq_1\right)\\
=\hspace{0.3em}
&\frac{1}{n} + \frac{1}{(n-i)n} \frac{1-\tfrac{k}{n-i+k}}{1-\tfrac{k}{n-i+k}\tfrac{i-k}{n-\ell-k}} \frac{n-\ell}{n-\ell-k}\\
=\hspace{0.3em}
&\frac{1}{n} + \frac{1}{(n-i)n} \frac{n-i}{(n-i+k)(n-\ell-k)-k(i-k)} (n-\ell)\\
=\hspace{0.3em}
&\frac{1}{n}\left(1 +  \frac{n-\ell}{(n-\ell)(n-i)-k\ell}\right)\\
=\hspace{0.3em}
&\frac{1}{n}\frac{(n-\ell)(n-i+1)-k\ell}{(n-\ell)(n-i)-k\ell}.
\end{align*}
Hence, conditionally on the event that the path $v_1,\dots,v_{k+1}$ is already constructed, the probability that the path is extended to $v_1,\dots,v_{n-\ell+1}$ is
\begin{equation*}
\prod_{i=k+1}^{n-\ell} \frac{1}{n}\frac{(n-\ell)(n-i+1)-k\ell}{(n-\ell)(n-i)-k\ell} = \frac{1}{n^{n-\ell-k-1}\ell}.
\end{equation*}

Now, it remains to compute the probability that the path is extended to $v_n$. We do this in a similar fashion: consider the Markov chain $(M_{2,i})_{i\ge 0}$ with three states, $A_2$, $B_2$ and $C_2$, and transition probabilities 
\begin{align*}
p_2 =\hspace{0.3em} &P_2(A_2,B_2) = \tfrac{n-\ell-k}{n-\ell-k+(n-i)},\\
q_2 =\hspace{0.3em} &P(B_2,A_2) = \tfrac{i-(n-\ell)}{\ell},\\ 
&P(A_2,C_2) = 1-p_2,\\ 
&P(B_2,C_2) = 1-q_2,\\ 
&P(C_2,C_2) = 1.
\end{align*}
By the same logic as above, conditionally on the event that the path $v_1,\dots,v_i$ is already constructed, the probability that $v_iv_{i+1}$ is the next edge is
\begin{align*}
&\dfrac{1}{n-k} + \dfrac{1}{(i-(n-\ell))(n-i)}\left(\frac{i-(n-\ell)}{n-k}\frac{1-p_2}{1-p_2q_2} + \frac{n-\ell-k}{n-k}\frac{q_2(1-p_2)}{1-p_2q_2}\right)\\
=\hspace{0.3em}
&\frac{1}{n-k} + \frac{1}{(i+\ell-n)(n-i)} \frac{1-p_2}{(n-k)(1-p_2q_2)} \left(i+\ell-n + (n-\ell-k)q_2\right)\\
=\hspace{0.3em}
&\frac{1}{n-k} + \frac{1}{(n-i)(n-k)} \frac{1-\tfrac{n-\ell-k}{n-\ell-k+(n-i)}}{1-\tfrac{n-\ell-k}{n-\ell-k+(n-i)}\tfrac{i-(n-\ell)}{\ell}} \frac{n-k}{\ell}\\
=\hspace{0.3em}
&\frac{1}{n-k} + \frac{1}{(n-i)} \frac{n-i}{(2n-\ell-k-i)\ell-(n-\ell-k)(i+\ell-n)}\\
=\hspace{0.3em}
&\frac{1}{n-k} + \frac{1}{(n-i)(n-k)}\\
=\hspace{0.3em}
&\frac{n-i+1}{(n-i)(n-k)}.
\end{align*}
Taking the product from $i = n-\ell+1$ to $n-1$ yields $\tfrac{\ell}{(n-k)^{\ell-1}}$, and in total the probability that the path $v_1,\dots,v_n$ is sampled as a uniform spanning tree of $H_n$ is
\begin{equation*}
\frac{1}{(n-\ell)^{k-1}(n-k)^{\ell-1} (n-\ell-k) n^{n-\ell-k-1}},
\end{equation*}
which means that there are $(n-\ell)^{k-1}(n-k)^{\ell-1} (n-\ell-k) n^{n-\ell-k-1}$ spanning trees of $H_n$. One may conclude since the number of spanning trees of $K_n$ is $n^{n-2}$ by Cayley's formula.
\end{proof}

A slightly stronger version of the next lemma essentially appears as Theorem~10 in~\cite{JLR}. We provide a short proof for the sake of completeness.

\begin{lemma}\label{lem perm}
Fix $\eps = \eps(n)\in (0,1)$, $k = k(n)$ and $s = \lfloor \eps n\rfloor$. Let $\nu$ be a permutation of $[n]$ chosen uniformly at random and let $U\subseteq [n]$ satisfying $|U|=k$. Then,
$$\Prob(|\nu([s])\cap U|\le \eps k/2)\le 2\exp\left(-\dfrac{\eps k}{36}\right).$$
\end{lemma}
\begin{proof}
Fix $p = 3\eps/4$. Sample $n$ iid random variables $(X_i)_{i=1}^n$ with distribution $\mathrm{Unif}([0,1])$ and order them as $X_{i_1} < \dots < X_{i_n}$ (this ordering is a.s.\ well defined as ties do not appear a.s.). Then, define the permutation $\nu': j\in [n]\mapsto i_j\in [n]$. Clearly $\nu'\sim \nu$.

Now, by Chernoff's inequality (Lemma~\ref{chernoff})
\begin{equation*}
    \Prob(|\{i: X_i\le p\}|\ge \eps n)\le \exp\left(-\dfrac{3\eps n/4}{27}\right) = \exp\left(-\dfrac{\eps n}{36}\right),
\end{equation*}
and moreover
\begin{equation*}
    \Prob(|\{i\in U: X_i\le p\}|\le \eps k/2)\le \exp\left(-\dfrac{3\eps k/4}{27}\right) = \exp\left(-\dfrac{\eps k}{36}\right).
\end{equation*}
Thus, we conclude that
\begin{equation*}
\Prob(|\nu'([s])\cap U|\le \eps k/2)\le \Prob(|\{i: X_i\le p\}|\ge \eps n \text{ or }|\{i\in U: X_i\le p\}|\le \eps k/2)\le 2\exp\left(-\dfrac{\eps k}{36}\right),
\end{equation*}
which proves the lemma.
\end{proof}

\begin{lemma}\label{lem bipart AB}
Fix any sufficiently small $\varepsilon > 0$ and integers $m = m(n)\ge \eps n$, $\ell = \ell(n)\in [\eps m, m]$ and $k = k(n)$ satisfying $k+m\le n$. Fix a bipartite graph $G = (U,V,F)$, where $F = E(G)$, satisfying $U,V\subseteq V_n$, $|U| = k$, $|V| = m$ and $\forall u\in U, \deg_G(u)\ge \ell$. Let $T\sim \mathrm{Unif}(\cT_n)$. Then, there is a constant $c = c(\eps) > 0$ such that
\begin{equation*}
\Prob(F\cap E(T) = \emptyset)\le \exp(-ck).
\end{equation*}
\end{lemma}
\begin{proof}
Our main tool is again the Aldous-Broder algorithm. Let us start a lazy random walk $(W_t)_{t\ge 0}$ on $K_n$ from a uniformly chosen vertex. Define $n_1 = \eps^3 n$ and the event $\cB = \{|\{W_t\}_{t=0}^{n_1-1}| > n_1/2\}$. For every sufficiently small $\eps$, an immediate application of Azuma's inequality for the martingale $(|\{W_t\}_{t=0}^i|-\mathbb E|\{W_t\}_{t=0}^i|)_{i=0}^{n_1-1}$ (which changes by at most 1 with every step) shows that 
$$\mathbb P(\overline{\cB})\le \exp\left(-\dfrac{(n_1/4)^2}{2n_1}\right) = \exp\left(-\dfrac{\eps^3 n}{32}\right),$$
where the inequality holds since $\mathbb E|\{W_t\}_{t=0}^{\eps^3 n-1}| = (1+o_{\eps,n}(1))\eps^3 n$ (reflecting the idea that the number of terms, visited by $(W_t)_{t=0}^{\eps^3 n-1}$ more than once, becomes negligible when $\eps\to 0$).

For every $i\in [k]$, we denote $\tau_i = \min\{j\ge 0: |\{W_t\}_{t=0}^j\cap U| = i\}$ and we define the event $\cC_i = \{\tau_i\le n_1\}$. Let also $k_1 = \eps^3 k/4$. By Lemma~\ref{lem perm} for the set of vertices, met by the random walk until time $n_1$, we have
\begin{equation*}
\Prob(\cC_1)\ge \dots\ge \Prob(\cC_{k_1})\ge 1 - \Prob(\overline{\cC_{k_1}}\mid \cB) + \Prob(\cB)\ge 1 - \exp(-\Omega_{\eps}(k)) - \exp(-\Omega_{\eps}(n)) = 1 -  \exp(-\Omega_{\eps}(k)).
\end{equation*}

\noindent
Now, for every $i\in [k]$, define the events
\begin{equation*}
    \cE_i = \{W_{\tau_i+1}\in N_V(W_{\tau_i})\} \text{ and } \cF_i = \{W_{\tau_i+1}\notin (W_t)_{t=0}^{\tau_i-1}\}.
\end{equation*}
One may easily show that $\{F\cap E(T) = \emptyset\}\subseteq \cap_{i=1}^k \overline{\cE_i\cap \cF_i}$, so in the remainder of the proof we will bound from above the probability of the event $\cap_{i=1}^k \overline{\cE_i\cap \cF_i}$.

Note that, on the one hand, for all $i\in [k]$, $\mathbb P(\overline{\cE_i})\le 1-\ell/n\le 1-\eps^2$ and $\cE_i$ is independent from $(W_t)_{t=0}^{\tau_i-1}$. On the other hand, for any event $\cD$ measurable with respect to the $\sigma$-algebra generated by $(W_t)_{t=0}^{\tau_i}$, $\mathbb P(\overline{\cF_i}\mid \cC_i\cap \cD)\le \eps^3$. We deduce that
\begin{align*}
&\Prob\left(\cap_{i=1}^{k} \overline{\cE_i\cap \cF_i}\right)\\
\le\hspace{0.3em} 
&\Prob\left(\cap_{i=1}^{k_1} \overline{\cE_i\cap \cF_i}\right)\\
\le\hspace{0.3em} 
&\Prob\left(\cap_{i=1}^{k_1} \left(\overline{\cE_i} \cup \left(\overline{\cF_i}\cap \cC_i\right)\right)\right) + \Prob\left(\cup_{i=1}^{k_1} \overline{\cC_i}\right)\\
\le\hspace{0.3em} 
&\prod_{i=1}^{k_1} \left(\Prob\left(\overline{\cE_i}\mid \cap_{j=1}^{i-1} \left(\overline{\cE_j} \cup \left(\overline{\cF_j}\cap \cC_j\right)\right)\right) + \Prob\left(\overline{\cF_i}\cap \cC_i\mid \cap_{j=1}^{i-1} \left(\overline{\cE_j} \cup \left(\overline{\cF_j}\cap \cC_j\right)\right)\right)\right) + k\exp(-\Omega_{\eps}(k))\\
\le\hspace{0.3em} 
&\prod_{i=1}^{k_1} \left(\Prob\left(\overline{\cE_i}\right) + \Prob\left(\overline{\cF_i}\mid \cC_i\cap \left(\cap_{j=1}^{i-1} \left(\overline{\cE_j} \cup \left(\overline{\cF_j}\cap \cC_j\right)\right)\right)\right)\right) + k\exp(-\Omega_{\eps}(k))\\
\le\hspace{0.3em} 
&\prod_{i=1}^{k_1} \left(1 - n^{-1}\ell + \eps^3\right) + k\exp(-\Omega_{\eps}(k))\\
\le\hspace{0.3em} 
&\left(1 - \eps^2 + \eps^3\right)^{k_1} + k\exp(-\Omega_{\eps}(k)) = \exp(-\Omega_{\eps}(k)),
\end{align*}
which concludes the proof of the lemma.
\end{proof}

\section{Proofs of the main results}\label{sec proofs}

\subsection{Proof of Theorem~\ref{thm min degree}}

We begin with a proof of the lower bound. It is an application of the differential equation method from the moment when there are $n^{1/2}\log n$ vertices of degree at most $k-1$ left to the moment when there are none.

\begin{proof}[Proof of Theorem~\ref{thm min degree}~\eqref{pt 1 thm deg}]

Fix $\eps\in (0,1)$ and any strategy $\sigma\in \mathcal S_n$. Denote by $s_{i,k}$ the number of vertices of degree at most $k-1$ in $G_i$ and define $i_k = i_k(\eps) = \min\{i\in \mathbb N: s_{i,k}\le \eps^{-1} n^{1/2}\}$. Thus, at every round $i\ge i_k$, 
\begin{equation}\label{eq bounded 1}
s_{i,k}-s_{i+1,k}\le 2,
\end{equation}
and by Lemma~\ref{lem clique AB},
\begin{align}
\mathbb P(s_{i,k} - s_{i+1,k} 
= 2\mid s_{i,k})
&\le\hspace{0.3em}  \mathbb P\left(E(T_{i+1})\cap \binom{V_{s_{i,k}}}{2}\neq \emptyset\;\bigg|\; s_{i,k}\right)\label{eq:ineq}\\
&=\hspace{0.3em} 1 - \left(1-\tfrac{s_{i,k}}{n}\right)^{s_{i,k}-1}\nonumber\\
&=\hspace{0.3em} 1 - \exp\left(-\tfrac{s_{i,k}^2}{n} + O\left(\tfrac{s_{i,k}^3}{n^2}+\tfrac{s_{i,k}}{n}\right)\right)\nonumber\\
&=\hspace{0.3em} 1 - \exp\left(-\tfrac{s_{i,k}^2}{n} + O(n^{-1/2})\right),\nonumber
\end{align}
In fact,~\eqref{eq:ineq} holds with equality for strategies that connect two vertices of degree $k-1$ whenever it is possible. Hence,
\begin{equation}\label{eq expectation}
\mathbb E[s_{i,k}-s_{i+1,k}\mid s_{i,k}]\le 2 - \exp\left(-\tfrac{s_{i,k}^2}{n} + O(n^{-1/2})\right) = 2 - \exp\left(-\tfrac{s_{i,k}^2}{n}\right) + O(n^{-1/2}).
\end{equation}
Now, we define
\[x = x(\eps, n) = \min\left\{\frac{i - i_k}{n^{1/2}}: \exists \sigma\in \cS_n: \text{ all but at most }\eps n^{1/2} \text{ vertices of }G_{\sigma, i} \text{ have degree at least }k\right\}.\]
Now, we apply Theorem~\ref{Thm:DEMethod} with scaling of $n^{1/2}$ instead of $n$,  
$a=1$, $Y_{1}(i) = s_{i,k}$, $t = \tfrac{i-i_k}{n^{1/2}}$, $\alpha(t) = \tfrac{s_{i,k}}{n^{1/2}}$ and 
\[F_1\left(\frac{i-i_k}{n^{1/2}}, \frac{s_{i,k}}{n^{1/2}}\right) = 2 - \exp\left(-\frac{s_{i,k}^2}{n}\right).\]
In fact, using~\eqref{eq bounded 1} (verifying the boundedness hypothesis) and~\eqref{eq expectation} (verifying the trend hypothesis), the unique solution of the differential equation 
\[y' = -(2 - \exp(-y^2)) \text{ with } y(0) = \eps^{-1}\]
whp dominates $\alpha$ in the sense that $\alpha(t)\le y(t) + o(1)$ for all $t > 0$. 

In fact, obtaining a lower bound for $x$ does not require the resolution of the above differential equation. 
Indeed, by using that $x = \min\{(i-i_k)n^{-1/2}: s_{i, k} - s_{i_k, k} = (\eps^{-1}-\eps)n^{1/2}\}$, Theorem~\ref{Thm:DEMethod} implies that whp $x$ may be bounded from below (up to lower order terms coming from the controlled approximation error) via the inequality
\begin{equation*}
\int_{t=0}^{x} 2 - \exp(-y^2) dy\ge \eps^{-1}-\eps-o(1),
\end{equation*}
which leads to 
\begin{equation*}
2x - \left(\tfrac{(2\pi)^{1/2}}{2} - O(\exp(-x^2))\right) \ge \eps^{-1}-\eps-o(1),
\end{equation*}
or equivalently $x \ge \tfrac{\eps^{-1}-\eps}{2} + \tfrac{(2\pi)^{1/2}}{4} + o(1)$. Letting $\eps \to 0$ shows that whp the property of minimum degree $k$ cannot be achieved by Builder  within less than 
\[\frac{kn - \eps^{-1} n^{1/2}}{2} + \frac{\eps^{-1} -\eps+(2\pi)^{1/2}/2}{2}n^{1/2} + o(n^{1/2}) = \frac{kn}{2}+\frac{(2\pi)^{1/2}}{4}n^{1/2}+o_{\eps, n}(n^{1/2})\]
steps, as desired.
\end{proof}




The strategy we choose to follow in the proof of Theorem~\ref{thm min degree}~\eqref{pt 2 thm deg} and~\eqref{pt 3 thm deg} is described by the following algorithm. The idea of the proof is similar to the one of~\eqref{pt 1 thm deg} but this time we have to take care of possible repetitions of edges.

\begin{algorithm}\label{algo degree k}
\emph{Input: an empty graph on $n$ vertices. Initiate $i\leftarrow 1$. While $i < k+1$, given a set of vertices of degrees $i-1$, $i$ and $i+1$, do:
\begin{enumerate}[(i)]
    \item\label{step 1} if there is an edge, included in the uniform spanning tree proposed at the current step, which connects two vertices of degree $i-1$ that were not adjacent before, choose such an edge uniformly at random and construct it;
    \item\label{step 2} else, if there is an edge, included in the uniform spanning tree proposed at the current step, which connects a vertex of degree $i-1$ and a vertex of degree $i$ that were not adjacent before, choose such an edge uniformly at random and construct it;
    \item\label{step 3} else, if there is still at least one vertex of degree $i-1$ and there is an edge, included in the uniform spanning tree proposed at the current step, which connects two vertices of degree $i$ that were not adjacent before, choose such an edge uniformly at random and construct it;
    \item\label{step 4} else, if there are no vertices of degree $i-1$ left, do $i\leftarrow i+1$;
    \item\label{step 5} else, exit the loop.
\end{enumerate}
Output: the obtained graph.}
\end{algorithm}

Theorem~\ref{thm min degree}~\eqref{pt 2 thm deg} may be easily deduced by the following more general lemma after a union bound and the fact that Algorithm~\ref{algo degree k} produces a graph with maximum degree $k+1$ in case step~\eqref{step 5} is never executed. For every $m\ge 1$, denote by $\cE_m$ the event ``the iteration of the loop that enforces $i\leftarrow m$ is attained, and at this point there are at most $\tfrac{(2\pi)^{1/2}}{4} n^{1/2} + \tfrac{n^{1/2}}{\log n}$ vertices of degree $m+1$''.

\begin{lemma}\label{lemma key}
Fix $k = k(n) = o(n^{1/2})$. Then, $\Prob(\cE_1) = 1$, and for every integer $m\in [2,k+1]$, 
$$\Prob(\cE_m\mid \cE_{m-1}) = 1 - o(n^{-1}).$$
\end{lemma}

\begin{proof}[Proof of Lemma~\ref{lemma key}]
For $m=1$ the statement is trivial since the step $i\leftarrow 1$ happens when the graph contains no edges. 

Fix $m\in [2,k+1]$ and condition on the event $\cE_{m-1}$ and on the (random) graph $G$, observed at the round when $i\leftarrow m$ is implemented. Below we abuse notation and, despite the additional edges $G$ receives at every round (so formally we are talking about the sequence of graphs $(G_i)_{i\ge 0}$), we sometimes denote by $G$ the graph throughout the process for convenience. For all $j\ge 0$, denote by $U_{i,j}\subseteq V_n$, or just by $U_j$, the set of vertices of degree $j$ in $G_i$. Also, fix $\omega = (k^{-2}n)^{1/6}$. Our proof requires a consideration of four different regimes: $|U_{m-1}|\ge s := \lfloor (\log n)^3 n^{1/2}\rfloor$, $|U_{m-1}|\in [\omega^{-1} n^{1/2}, s]$, $|U_{m-1}|\in [r := \lfloor (\log n)^2\rfloor, \omega^{-1} n^{1/2}]$ and $|U_{m-1}|\in [1, r]$.

\paragraph{Regime 1.} We will show that as long as $|U_{m-1}|\ge s$, at every round, step~\eqref{step 1} of Algorithm~\ref{algo degree k} is executed with probability $1-o(n^{-2})$. Indeed, suppose that $U_{m-1}\subseteq V_n$ contains at least $s$ vertices. Colour the vertices in $U_{m-1}$ into $r$ colours uniformly at random and independently, and let $W_1,\dots,W_r$ be the respective colour classes. Then, an application of Chernoff's bound (Lemma~\ref{chernoff}) shows that for every $j\in [r]$, the events 
$$\cA_j = \{|W_j|\ge (2r)^{-1}s\} \text{ and } \cB_j = \{\text{the maximum degree of }G[W_j]\text{ is less than }\tfrac{2}{3}|W_j|-1\}$$ 
both hold with probability $1 - o(r^{-1})$.

For $T\sim \mathrm{Unif}(\mathcal T_n)$ and all $j\in [r]$, define the event $\cC_j = \left\{E(T)\cap \tbinom{W_j}{2}\neq \emptyset\right\}$. Then, by Lemma~\ref{lem clique AB}, for all $j\in [r]$,
$$\mathbb P(\cC_j)\ge \mathbb P(\cC_j\mid \cA_j)\mathbb P(\cA_j)\ge \left(1 - \left(1-\tfrac{s}{2rn}\right)^{(2r)^{-1} s-1}\right) (1-o(r^{-1})) = 1-o(r^{-1}).$$
\noindent
Now, conditionally on $\cC_j$, sample an edge $f_j\in E(T)\cap \tbinom{W_j}{2}$ uniformly at random. Since by symmetry $f_j\sim \mathrm{Unif}\left(\tbinom{W_j}{2}\right)$, conditionally on $\cB_j\cap \cC_j$ the edge $f_j$ has not been constructed yet with probability at least $1/3$. Thus, the probability that step~\eqref{step 1} of Algorithm~\ref{algo degree k} is not executed is bounded from above by
\begin{equation*}
    \left(1-\tfrac{1}{3}\mathbb P(\cap_{j=1}^r (\cB_j\cap \cC_j))\right)^r \le \left(1-\tfrac{1}{3}(1-o(1))\right)^r = o(n^{-2}),
\end{equation*}
which is sufficient to conclude in this case.\\

Before continuing with the proof, let us point out that the second regime is key for understanding where the expression of the second term of $\tau(\sigma, \cP_{\min, k})$ comes from. Unlike the proof of the lower bound where we were allowed to ignore edges of $G[U_{k-1}]$, here we need to take more care to ensure that, roughly speaking, the process ``behaves similarly'' with or without the edges of $G[U_{m-1}]$. It is at this point that the assumption of $k = o(n^{1/2})$ is used.

\paragraph{Regime 2.} Denote by $s_{i,m-1}$ the size of $U_{m-1}$ at round $i$. As in the proof of the lower bound, let $x(n) n^{1/2}$ be the number of rounds needed for Algorithm~\ref{algo degree k} to go from $s_{i,m-1} = s$ to $s_{i,m-1} = \omega k$ vertices of degree $m-1$ in $G$. Also, set $i_{m-1} = \min\{i: s_{i,m-1}\le s\}$. Observe that Lemma~\ref{lem bipart AB} for the complementary graph of $G[U_{m-1}, V_n\setminus (U_{m-1}\cup U_{m+1})]$ ensures that, at any round before reaching $s_{i, m-1} = \omega k$, Algorithm~\ref{algo degree k} executes either step~\eqref{step 1} or step~\eqref{step 2} with probability at least $1 - \exp(-\Omega(\omega k)) = 1 - o(n^{-2})$.

Before we compute the expected one-round changes of $(s_{i,m-1})_{i\ge i_{m-1}}$, let $(U_{m-1}^{(1)}, U_{m-1}^{(2)})$ be an arbitrary partition of the set $U_{m-1}$. Then, on the one hand,
\begin{align}
&\mathbb P(s_{i,m-1} - s_{i+1,m-1} 
= 2\mid s_{i,m-1})\nonumber\\
=\hspace{0.3em}  
&\mathbb P\left(E(T_{i+1})\cap (\tbinom{U_{i,m-1}}{2}\setminus E(G_i))\neq \emptyset\mid s_{i,m-1}\right)\nonumber\\
\ge\hspace{0.3em}  
&\left(1 - \mathbb P(E(T_{i+1})\cap \tbinom{U_{i,m-1}^{(1)}}{2} = \emptyset\mid s_{i,m-1}) \mathbb P(E(T_{i+1})\cap \tbinom{U_{i,m-1}^{(2)}}{2} = \emptyset\mid s_{i,m-1})\right) \left(1-\left(\tfrac{k}{s_{i,m-1}/2-2}\right)^2\right)\label{eq explain 1}\\
\ge\hspace{0.3em}
&\left(1 - \left(1-\tfrac{s_{i,m-1}/2-2}{2n}\right)^{2(s_{i,m-1}/2-3)}\right) \left(1-\tfrac{16k^2}{s_{i,m-1}^2}\right)\label{eq explain 2}\\
=\hspace{0.3em} 
&\left(1 - \exp\left(-\tfrac{s_{i,m-1}^2}{4n} + O\left(\tfrac{s_{i,k}^3}{n^2}+\tfrac{s_{i,k}}{n}\right)\right)\right) \left(1-\tfrac{16k^2}{n}\left(\tfrac{s_{i,m-1}^2}{4n}\right)^{-1}\right)\nonumber\\
\ge\hspace{0.3em} 
&\left(1 - \exp\left(-\tfrac{s_{i,m-1}^2}{4n} + O((\log n)^9 n^{-1/2})\right)\right) \left(1-16\omega^{-6} \left(\tfrac{s_{i,m-1}^2}{4n}\right)^{-1}\right)\nonumber\\
=\hspace{0.3em} 
&1 - 16\omega^{-6} \left(\tfrac{s_{i,m-1}^2}{4n}\right)^{-1} - \left(1-16\omega^{-6} \left(\tfrac{s_{i,m-1}^2}{4n}\right)^{-1}\right)\exp\left(-\tfrac{s_{i,m-1}^2}{4n} + O((\log n)^9 n^{-1/2})\right),\label{eq add factor 1}
\end{align}
where~\eqref{eq explain 1} and~\eqref{eq explain 2} are due to the fact that the the proportion of edges to the total number of vertex pairs in both $G[U_{m-1}^{(1)}]$ and $G[U_{m-1}^{(2)}]$ is always at most $\tfrac{k}{\lfloor s_{i,m-1}/2\rfloor - 1}\le \tfrac{k}{s_{i,m-1}/2-2}$. On the other hand, a similar computation without dividing $U_{m-1}$ leads to
\begin{align}
&\mathbb P(s_{i,m-1} - s_{i+1,m-1} 
= 2\mid s_{i,m-1})\nonumber\\
=\hspace{0.3em}  
&\mathbb P\left(E(T_{i+1})\cap (\tbinom{U_{i,m-1}}{2}\setminus E(G_i))\neq \emptyset\mid s_{i,m-1}\right)\nonumber\\
\ge\hspace{0.3em}  
&\left(1 - \mathbb P\left(E(T_{i+1})\cap \tbinom{U_{i,m-1}}{2}\neq \emptyset\mid s_{i,m-1}\right)\right) \left(1-\tfrac{k}{s_{i,m-1}-1}\right)\label{eq explain}\\
\ge\hspace{0.3em}
&\left(1 - \left(1-\tfrac{s_{i,m-1}}{n}\right)^{s_{i,m-1}-1}\right) \left(1-\tfrac{2k}{s_{i,m-1}}\right)\nonumber\\
=\hspace{0.3em} 
&\left(1 - \exp\left(-\tfrac{s_{i,m-1}^2}{n} + O\left(\tfrac{s_{i,k}^3}{n^2}+\tfrac{s_{i,k}}{n}\right)\right)\right) \left(1-\tfrac{2k}{n^{1/2}}\left(\tfrac{s_{i,m-1}^2}{n}\right)^{-1/2}\right)\nonumber\\
\ge\hspace{0.3em} 
&\left(1 - \exp\left(-\tfrac{s_{i,m-1}^2}{n} + O((\log n)^9 n^{-1/2})\right)\right) \left(1-2\omega^{-3} \left(\tfrac{s_{i,m-1}^2}{n}\right)^{-1/2}\right)\nonumber\\
=\hspace{0.3em} 
&1 - 2\omega^{-3} \left(\tfrac{s_{i,m-1}^2}{n}\right)^{-1/2} - \left(1-2\omega^{-3} \left(\tfrac{s_{i,m-1}^2}{n}\right)^{-1/2}\right)\exp\left(-\tfrac{s_{i,m-1}^2}{n} + O((\log n)^9 n^{-1/2})\right).\label{eq add factor}
\end{align}
Also, by Lemma~\ref{lem bipart AB},
\begin{align*}
\mathbb P(s_{i,m-1} - s_{i+1,m-1} 
= 1\mid s_{i,m-1}) = 1 - \mathbb P(s_{i,m-1} - s_{i+1,m-1} 
= 2\mid s_{i,m-1}) - \exp(-\Omega(\omega k)),
\end{align*}
and hence, using both~\eqref{eq add factor 1} when $s_{m-1}\ge \omega n^{1/2}$ and~\eqref{eq add factor} when $s_{m-1} < \omega n^{1/2}$, we obtain
\begin{align}
& \mathbb E[s_{i,m-1} - s_{i+1,m-1}\mid s_{i,m-1}]\nonumber\\
=\hspace{0.3em} 
& 1 + \mathbb P(s_{i,m-1} - s_{i+1,m-1} 
= 2\mid s_{i,m-1}) - \exp(-\Omega(\omega k))\nonumber\\
\ge\hspace{0.3em} 
& 2 - \left(\exp\left(-\tfrac{s_{i,m-1}^2}{4n} + (\log n)^9 n^{-1/2}\right) + O\left(\omega^{-6} \left(\tfrac{s_{m-1}^2}{n}\right)^{-1}\right)\right) \mathds{1}_{s_{m-1}\; \ge\; \omega n^{1/2}}\nonumber\\
& \hspace{0.5em} - \left(\exp\left(-\tfrac{s_{i,m-1}^2}{n} + (\log n)^9 n^{-1/2}\right) + O\left(\omega^{-3} \left(\tfrac{s_{i,m-1}^2}{n}\right)^{-1/2}\right)\right) \mathds{1}_{s_{i,m-1}\; <\; \omega n^{1/2}}.\label{eq expectation 2}
\end{align}

Again, renormalising time by setting $t = \tfrac{i-i_{m-1}}{n^{1/2}}$ and writing $s_{i,m-1} = \alpha n^{1/2} = \alpha(t) n^{1/2}$, one may deduce (as in the proof of Theorem~\ref{thm min degree}~\eqref{pt 1 thm deg} but with equality in the other direction) that
\begin{equation*}
\int_{\alpha=\omega^{-1}}^{x} \eqref{eq expectation 2}\; d\alpha\le \omega,
\end{equation*}
which can be rewritten as
\begin{equation*}
2x - \int_{\omega}^{x} \left(\exp\left(-\tfrac{\alpha^2}{4} + n^{-1/2+o(1)}\right) + O\left(\tfrac{1}{\omega^6\alpha^2}\right)\right) d\alpha - \int_{\omega^{-1}}^{\omega} \left(\exp\left(-\alpha^2 + n^{-1/2+o(1)}\right) + O\left(\tfrac{1}{\omega^3\alpha}\right)\right) d\alpha \le \omega.
\end{equation*}
Taking into account that the first integral is $o(1)$ and the second integral is given by 
$$(1+o(1)) \int_{0}^{\infty} \exp(-\alpha^2) d\alpha + O(\tfrac{\log \omega}{\omega^3}) \le \tfrac{(2\pi)^{1/2}}{2}+o(1),$$ we conclude that, with probability $1 - o(n^{-1})$,
\begin{equation*}
x \le \tfrac{1}{2} n^{-1/2} s + \tfrac{(2\pi)^{1/2}}{2}+o(1).
\end{equation*}

Consequently, 
when $s_{m-1}$ becomes at most $\omega k$, the number of vertices of degree $m+1$ is at most $\tfrac{(2\pi)^{1/2}}{2} n^{1/2} + o(n^{1/2})$ with probability $1 - o(n^{-1})$.

\paragraph{Regime 3.} In this regime $|U_{m+1}| = o(n^{1/2})$. Then, by consecutive applications of Lemma~\ref{lem bipart AB} for the complementary graph of $G[U_{m-1}, V_n\setminus (U_{m-1}\cup U_{m+1})]$ (note that the second part has size $n-o(n)$), at every round either step~\eqref{step 1} or step~\eqref{step 2} of Algorithm~\ref{algo degree k} is executed with probability $1-o(n^{-2})$. Thus, by a union bound over $\omega^{-1} n^{1/2} - r = O(n)$ rounds, there are at most $\tfrac{(2\pi)^{1/2}}{2} n^{1/2} + o(n^{1/2}) + \omega^{-1} n^{1/2} = \tfrac{(2\pi)^{1/2}}{2} n^{1/2} + o(n^{1/2})$ vertices of degree $m+1$ with probability $1-o(n^{-1})$ until the end of this regime.

\paragraph{Regime 4.} Finally, when $s_{m-1} \le r$, since $m=o(n^{1/2})$, the probability that a vertex in $U_{m-1}$ is connected by an edge (which is not yet in $G$) of the uniform spanning tree of $K_n$, sampled at the current step, to a vertex in $V_n\setminus (U_{m-1}\cup U_{m+1})$ is $1 - o(n^{-1/2}+n^{-1}(\log n)^2) = 1-o(n^{-1/2})$. Thus, $10 r = o(n^{1/2})$ independent trials of Algorithm~\ref{algo degree k} succeed to empty the set $U_{m-1}$ with probability $1 - r (n^{-1/2})^{10} = 1-o(n^{-1})$, which concludes the proof of the lemma by a union bound over all four regimes.
\end{proof}

To show Theorem~\ref{thm min degree}~\eqref{pt 3 thm deg}, we state a modification of Lemma~\ref{lemma key} and indicate the differences with the above proof. For every $m\ge 1$, denote by $\cE_m'$ the event ``the iteration of the loop that enforces $i\leftarrow m$ is attained, and at this point there are at most $O(\max\{m, n^{1/2}\})$ vertices of degree $m+1$''.

\begin{lemma}\label{lemma key mod}
Fix $k = k(n)$ satisfying $k = \Omega(n^{1/2})$ and $k = o(n)$. Then, $\Prob(\cE_1') = 1$, and for every integer $m\in [2,k+1]$, 
\begin{equation*}
\Prob(\cE_m'\mid \cE_{m-1}') = 1 - o(n^{-1}).
\end{equation*}
\end{lemma}

\begin{proof}[Sketch of proof of Theorem~\ref{thm min degree}~\eqref{pt 3 thm deg}]
As before, clearly $\Prob(\cE_1') = 1$. Consider the following three regimes:
\begin{align*}
&|U_{m-1}|\in [\min\{(\log n)^3 k, n\}, n],\\ 
&|U_{m-1}|\in [4k, \min\{(\log n)^3 k, n\}],\\
&|U_{m-1}|\in [1,4k].
\end{align*}
The first regime is treated in verbatim the same way as Regime 1 above, and in particular, at every round while $|U_{m-1}|\ge (\log n)^3 k$, step~\eqref{step 1} of Algorithm~\ref{algo degree k} is executed with probability $1-o(n^{-2})$. Moreover, no vertices of degree $m+1$ appear while $|U_{m-1}|\ge \min\{(\log n)^3 k, n\}$.

The second regime is treated in precisely the same way as the first part of Regime 2 above with the modification that $s_{m-1} = \alpha \min\{m, n^{1/2}\}$. Indeed, since here we do not aim for an exact constant as before, the key factor for us is that the term $\alpha^{-2}$ in the expression of the error bound is integrable with respect to $\alpha$ on any interval bounded away from 0. 

The third regime is treated in verbatim the same way as Regime 3 above. Note that the assumption that $k = o(n)$ is used to allow the application of Lemma~\ref{lem bipart AB}.
\end{proof}

\subsection{\texorpdfstring{Proof of Theorems~\ref{thm k-conn}}{}}\label{sec k-conn}

From this point on, we often tacitly assume different arithmetic conditions (mostly concerning divisibility): in fact, these are never important for our argument. Also, lower and upper integer parts are spared when rounding is not important for the argument.

\begin{lemma}\label{lem bipartite}
Fix any function $\omega = \omega(n)$ growing to infinity as $n\to \infty$ arbitrarily slowly, and fix $s\in [\omega^2 n^{3/4}, n/2]$. Fix any two disjoint subsets $A,B\subseteq [n]$, both of size $s$. Then, there is a strategy $\sigma\in \cS_n$ that constructs a perfect matching between the two sets in time at most $s + O(\omega n^{3/4})$ whp.
\end{lemma}
\begin{proof}
Denote $m = \omega n^{1/2}$. We divide the construction into three stages.

The first stage consists of the first $r = s - m/2$ rounds. At each of these rounds, construct an edge between a vertex in $A$ and a vertex in $B$ that are still unmatched, if this is possible, else ignore the round and proceed further. By Lemma~\ref{lem 2.8} applied at round $i\le r$ with $\ell = k = s-i+1$, a round is ignored with probability at most $\exp\left(-(1+o(1))\tfrac{k^2+\ell^2}{n}\right) = \exp\left(-\Omega\left(\tfrac{(s-i)^2}{n}\right)\right)$. By Chernoff's inequality for the random variables $(X_i)_{i=1}^{r} = (\mathds{1}_{\text{round }i\text{ is ignored}})_{i=1}^r$ and the fact that
\begin{equation*}
\mathbb E\left[\sum_{i=1}^{r} X_i\right]\le \sum_{j=m/2}^{\infty} \exp\left(-\Omega\left(\frac{j^2}{n}\right)\right) = O(\exp(-\Omega(\omega^2)) m),
\end{equation*}
whp no more than $m/2$ rounds are ignored during the first $r$ rounds. We condition on this event in the sequel.




Define $A'$ as the set of matched vertices in $A$ after the first stage and $A'' = A\setminus A'$. Based on the first stage, we know that $|A''|\le m$, and we may safely assume that $|A''| = m$ by ignoring certain edges in the matching if necessary. 
Thus, let $a''_1, \dots , a''_{m}$ be an arbitrary ordering of the vertices of $A''$. Define $B'$ and $B''$ analogously. Our goal at the second stage will be to construct $m$ vertex-disjoint stars with centers $a''_1,\dots,a''_m$, each having $n^{1/4}$ leaves in $B'$. 
To do this, at each round $i\in [r+1, r+5mn^{1/4}]$, we choose uniformly at random and independently from previous rounds an edge in $E(T_i)\cap (A''\times B')$, if there are any, and construct it, if it has not been added yet. In fact, by Lemma~\ref{lem 2.8} and a union bound there exists such an edge at any of the $5mn^{1/4}$ rounds whp since $5mn^{1/4} \left(1-\tfrac{m}{n}\right)^{s-m} \left(1-\tfrac{s-m}{n}\right)^{m} = \exp(-\Omega(sm/n)) = \exp(-\Omega(n^{1/4}))$. We condition on this event. 
Moreover, the probability that we attempt to construct the same edge twice is bounded from above by $\tbinom{5mn^{1/4}}{2} (sm)^{-2} = o(n^{-1})$, so an immediate second moment computation shows that there are no more than $mn^{1/4}$ such edges whp. 
In particular, after $5mn^{1/4}$ steps the obtained graph dominates whp the binomial random bipartite graph between $A''$ and $B'$ where every edge appears independently with probability $\tfrac{3n^{1/4}}{s}$. By standard concentration arguments (see e.g. Chapter 3 in~\cite{FK16}) 
we obtain that whp each vertex $a''_i\in A''$ is adjacent to at least $2 n^{1/4}$ vertices in $B'$, and at least $n^{1/4}$ of its neighbours are not incident to any other vertex in $A''$ (indeed, the probability that a fixed neighbour $b'\in B'$ of $a''_i$ has another neighbour in $A''$ is $O(m n^{1/4}/s) = o(1)$). 
This ensures that, at the end of this stage, one may find whp a subgraph, consisting of $m$ disjoint stars centered at $(a''_i)_{i=1}^m$, each with $n^{1/4}$ leaves in $B'$. We condition on this event.

\begin{figure}
\centering
\begin{tikzpicture}[scale=0.8,line cap=round,line join=round,x=1cm,y=1cm]
\definecolor{wewewe}{rgb}{0.7,0.7,0.7}
\clip(-7,-6.081930832044296) rectangle (10.3,5.172752421323636);
\draw [line width=0.8pt] (-5,4)-- (-6,3);

\draw [line width=0.1pt] (-40,-1.5)-- (10.5,-1.5);
\draw [line width=0.1pt] (3.5,-1)-- (3.5,6);
\draw [line width=0.1pt] (3.5,-10)-- (3.5,-2);
\draw [rotate around={0:(-0.5,2)},line width=0.8pt] (-0.5,2) ellipse (3.6138621999185303cm and 0.9cm);
\draw [rotate around={0:(-0.5,0)},line width=0.8pt] (-0.5,0) ellipse (3.613862199918533cm and 0.9cm);
\draw [rotate around={0:(8,2)},line width=0.8pt] (8,2) ellipse (2.1540659228538095cm and 0.8cm);
\draw [rotate around={0:(8,0)},line width=0.8pt] (8,0) ellipse (2.154065922853807cm and 0.8cm);

\draw [rotate around={0:(-0.5,2-5)},line width=0.8pt] (-0.5,2-5) ellipse (3.6138621999185303cm and 0.9cm);
\draw [rotate around={0:(-0.5,0-5)},line width=0.8pt] (-0.5,0-5) ellipse (3.613862199918533cm and 0.9cm);
\draw [rotate around={0:(8,2-5)},line width=0.8pt] (8,2-5) ellipse (2.1540659228538095cm and 0.8cm);
\draw [rotate around={0:(8,0-5)},line width=0.8pt] (8,0-5) ellipse (2.154065922853807cm and 0.8cm);

\draw [line width=0.8pt] (-5,4)-- (-5,3);
\draw [line width=0.8pt] (-5,4)-- (-4,2);
\draw [line width=0.8pt] (-4,2)-- (-4,4);
\draw [line width=0.8pt] (-4,2)-- (-4,0);
\draw [line width=0.8pt] (-2,4)-- (-3,3);
\draw [line width=0.8pt] (-2,4)-- (-2,3);
\draw [line width=0.8pt] (-2,4)-- (-1,2);
\draw [line width=0.8pt] (-1,2)-- (-1,4);
\draw [line width=0.8pt] (-1,2)-- (-1,0);
\draw [line width=0.8pt] (2,4)-- (1,3);
\draw [line width=0.8pt] (2,4)-- (2,3);
\draw [line width=0.8pt] (2,4)-- (3,2);
\draw [line width=0.8pt] (3,2)-- (3,4);
\draw [line width=0.8pt] (3,2)-- (3,0);
\draw [line width=0.8pt] (5,4)-- (4,3);
\draw [line width=0.8pt] (5,4)-- (5,3);
\draw [line width=0.8pt] (5,4)-- (6,2);
\draw [line width=0.8pt] (6,2)-- (6,4);
\draw [line width=0.8pt] (9,4)-- (8,3);
\draw [line width=0.8pt] (9,4)-- (9,3);
\draw [line width=0.8pt] (9,4)-- (10,2);
\draw [line width=0.8pt] (10,2)-- (10,4);
\draw [line width=0.8pt] (6,-5)-- (3,-3);
\draw [line width=0.8pt] (6,-5)-- (-1,-3);
\draw [line width=0.8pt,dash pattern=on 0.8pt off 4pt] (6,-5)-- (1,-3);
\draw [line width=0.8pt,dash pattern=on 0.8pt off 4pt] (6,-3)-- (1,-5);
\draw [line width=0.8pt] (6,-3)-- (-1,-5);
\draw [line width=0.8pt] (6,-3)-- (-4,-5);
\draw [line width=0.8pt,color=wewewe] (1,-3)-- (1,-5);
\begin{scriptsize}
\draw [fill=black] (-5,4) circle (1.5pt);
\draw [fill=black] (-6,3) circle (1.5pt);
\draw [fill=black] (-5,3) circle (1.5pt);
\draw [fill=black] (-4,2) circle (1.5pt);
\draw [fill=black] (-4,4) circle (1.5pt);
\draw [fill=black] (-4,0) circle (1.5pt);
\draw [fill=black] (-2,4) circle (1.5pt);
\draw [fill=black] (-3,3) circle (1.5pt);
\draw [fill=black] (-2,3) circle (1.5pt);
\draw [fill=black] (-1,2) circle (1.5pt);
\draw [fill=black] (-1,4) circle (1.5pt);
\draw [fill=black] (-1,0) circle (1.5pt);
\draw [fill=black] (2,4) circle (1.5pt);
\draw [fill=black] (1,3) circle (1.5pt);
\draw [fill=black] (2,3) circle (1.5pt);
\draw [fill=black] (3,2) circle (1.5pt);
\draw [fill=black] (3,4) circle (1.5pt);
\draw [fill=black] (3,0) circle (1.5pt);
\draw [fill=black] (5,4) circle (1.5pt);
\draw [fill=black] (4,3) circle (1.5pt);
\draw [fill=black] (5,3) circle (1.5pt);
\draw [fill=black] (6,2) circle (1.5pt);
\draw [fill=black] (6,4) circle (1.5pt);
\draw [fill=black] (9,4) circle (1.5pt);
\draw [fill=black] (8,3) circle (1.5pt);
\draw [fill=black] (9,3) circle (1.5pt);
\draw [fill=black] (10,2) circle (1.5pt);
\draw [fill=black] (10,4) circle (1.5pt);
\draw [fill=black] (6,0) circle (1.5pt);
\draw [fill=black] (10,0) circle (1.5pt);
\draw [fill=black] (7,3) circle (0.5pt);
\draw [fill=black] (0,3) circle (0.5pt);
\draw [fill=black] (0.2,3) circle (0.5pt);
\draw [fill=black] (-0.2,3) circle (0.5pt);
\draw [fill=black] (7.2,3) circle (0.5pt);
\draw [fill=black] (6.8,3) circle (0.5pt);
\draw [fill=black] (6,-5) circle (1.5pt);
\draw [fill=black] (3,-3) circle (1.5pt);
\draw [fill=black] (-1,-3) circle (1.5pt);
\draw [fill=black] (6,-3) circle (1.5pt);
\draw [fill=black] (1,-3) circle (1.5pt);
\draw [fill=black] (1,-5) circle (1.5pt);
\draw [fill=black] (-1,-5) circle (1.5pt);
\draw [fill=black] (-4,-5) circle (1.5pt);

\draw (-4.5,2) node {\Large{$A'$}};
\draw (-4.5,0) node {\Large{$B'$}};
\draw (5.5,2) node {\Large{$A''$}};
\draw (5.5,0) node {\Large{$B''$}};
\end{scriptsize}
\end{tikzpicture}
\caption{On the top: the moment when the matching between the vertices of $A'$ and $B'$ is entirely constructed but the bipartite graph $(A'', B'')$ is still empty. On the bottom: the two dotted edges and the grey edge form an augmenting path, which allows to increase the matching between $A$ and $B$.}
\label{fig 1}
\end{figure}
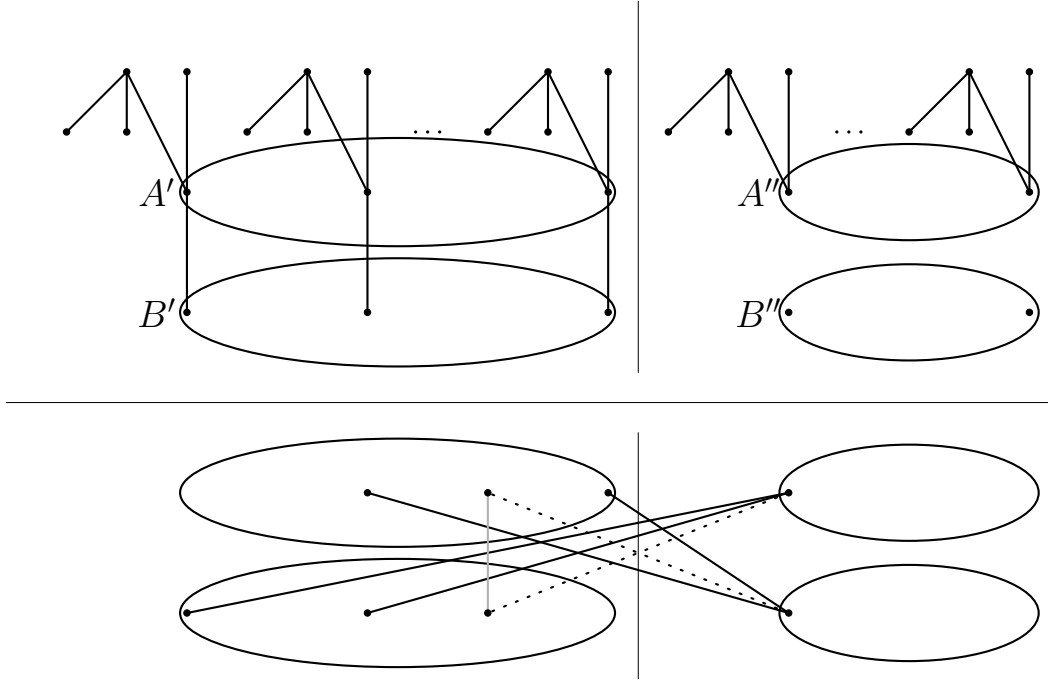

Now, for all $i\in [m]$, denote by $A_i$ the set of second neighbours of $a''_i$ and let $b''_1, \dots , b''_{m}$ be an ordering of the vertices of $B''$. At the third stage, our goal will be to (greedily) match the vertices of $B''$ with the sets $(A_i)_{i\in [m]}$ in the sense that, in the end, there must be a permutation $\pi$ of $[m]$ such that, for every $i\in [m]$, $b''_i$ is matched to a vertex of $A_{\pi(i)}$. Then, in the end of this third stage we will have at our disposal $m$ augmenting paths of length three, which will ensure the existence of a perfect matching (see the bottom part of Figure~\ref{fig 1}). For all $i\in [m]$, denote by $Z_i$ the number of rounds which, given that there are exactly $i$ unmatched vertices left in $B''$, are needed to connect an arbitrary unmatched vertex in $B''$ to an arbitrary of the $i$ corresponding sets among $(A_i)_{i=1}^m$. We will show that
\begin{equation*}
\sum_{i=1}^{m} \mathbb E[Z_i] = O(n^{3/4}),
\end{equation*}
which is sufficient to conclude by Markov's inequality for $\sum_{i=1}^m Z_i$. Indeed, fix $i\in [m]$ and suppose that there are $i$ still unmatched vertices in $B''$. By Lemma~\ref{lem 2.8} with $k = i n^{1/4}$ and $\ell = i$, the probability that none of the $i$ vertices is matched at the current step is
\begin{equation*}
\left(1 - \frac{i}{n}\right)^{i n^{1/4}} \left(1 - \frac{i n^{1/4}}{n}\right)^{i} \left(1 - \frac{i+ in^{1/4}}{n}\right) \le 1 - \frac{i^2}{n^{3/4}},
\end{equation*}
and hence $Z_i$ is dominated by a geometric random variable with parameter $i^2 n^{-3/4}$, which has expectation $n^{3/4} i^{-2}$. This is sufficient to conclude since the sum of $(i^{-2})_{i\ge 1}$ is finite.
\end{proof}

\begin{proof}[Proof of Theorem~\ref{thm k-conn}]
We define the strategy $\sigma$ as follows. To begin, iterate $k$ times the strategy, used in the proof of Lemma~\ref{lem bipartite}, to construct $k$ independent perfect matchings of $K_n$, and consider the union of these. Since some edges might coincide, we may obtain a graph with minimum degree less than $k$. However, we show that whp no vertex set of size between $2$ and $n/2$ of this graph can be disconnected from the rest by deleting at most $k-1$ other vertices. Finally, we add a few additional edges to ensure that the final graph has minimum degree $k$.

We proceed to the formal treatment of the problem. Fix $\omega = \omega(n) \to \infty$ as $n\to \infty$. By Lemma~\ref{lem bipartite} one may construct $k$ independent perfect matchings of $K_n$ in $\tfrac{kn}{2} + k \omega n^{3/4}$ rounds whp. Let $H_k$ be the graph obtained as a union of all $k$ perfect matchings. Note that another way to construct $H_k$ is the following: to every vertex in $V_n$ attach $k$ half-edges in $k$ different colours $(c_i)_{i=1}^k$. Then, match all half-edges uniformly at random so that every edge is matched with an edge of the same colour, and then identify repeated edges. (Note that every matching in a single color is uniformly random as its construction does not make use of the vertex labels, which may be distributed later.) Note that in the first moment computation below, repeated edges are always counted with multiplicities.

Let $U\subseteq V_n$ satisfy $s_U = |U|\in [2,n/2]$ and $W\subseteq V_n\setminus U$ satisfy $|W|=k-1$. We will show that whp $W$ does not disconnect $U$ from $V_n\setminus (U\cup W)$ for any choice of $(U,W)$. Having fixed $(U,W)$, a very simple upper bound on the probability of the above event, obtained by consecutively connecting monochromatic half-edges and analysing the probability at each step, is $\left(\tfrac{s_U+k-1}{n-1}\right)^{s_U k/2}$. Then, a union bound over all choices of $(U,W)$ leads to an upper bound of

\begin{equation}\label{eq union bd}
\tbinom{n}{s_U,k-1}\left(\tfrac{s_U+k-1}{n-1}\right)^{s_U k/2}\le \left(\tfrac{n}{s_U}\right)^{s_U} \left(\tfrac{n}{k-1}\right)^{k-1} \exp\left(s_U+k-1\right)\left(\tfrac{s_U+k-1}{n-1}\right)^{s_U k/2}.
\end{equation}

\noindent
Define $s_{\min} = s_{\min}(k)$ as $s_{\min}(3) = 7$, $s_{\min}(4) = 5$, $s_{\min}(5) = s_{\min}(6) = 4$ and, for all $k\ge 7$, $s_{\min}(k) = 3$. We first assume that $k\ge 3$, $s_U\ge s_{\min}$ and $s_U\le \eps n$ for a sufficiently small $\eps > 0$. Then,
\begin{align*}
\exp\left(s_U+k-1\right)\le \exp(s_Uk)\le \left(\tfrac{s_U+k-1}{n-1}\right)^{s_U k/\log(\eps^{-1})}
\end{align*}
as well as
\begin{align}
    \left(\tfrac{n}{s_U}\right)^{s_U}\left(\tfrac{n}{k}\right)^{k} \left(\tfrac{s_U+k-1}{n-1}\right)^{s_U k/2.01}
    &\le\hspace{0.3em} \left(1+\tfrac{k-1}{s_U}\right)^{s_U} \left(1+\tfrac{s_U-1}{k}\right)^{k} \left(\tfrac{s_U+k-1}{n-1}\right)^{s_U k/2.01 - s_U - k}\nonumber\\ 
    &\le\hspace{0.3em} \exp(s_U+k-2) \left(\tfrac{s_U+k-2}{n-1}\right)^{s_U k/2.01 - s_U - k}\nonumber\\
    &\le\hspace{0.3em} 
    \left(\tfrac{s_U+k-2}{n-1}\right)^{s_U k/2.01 - s_U - k - O\left(\tfrac{1}{\log \eps^{-1}}\right)(s_U+k)},\label{eq union bound 1}
\end{align}
where the last inequality is due to the fact that $\tfrac{n}{s_U}\ge (1+o(1))\eps^{-1}$. Now, for every sufficiently small $\eps > 0$ by assumption
\begin{align*}
\frac{s_U k}{2.01} - s_U - k - O\left(\tfrac{1}{\log \eps^{-1}}\right) (s_U+k)
&\ge\; \frac{s_U k - 2.03 s_U - 2.03 k}{2.01} + \dfrac{s_U+k}{201}\\
&=\; \frac{(s_U-2.03)(k - 2.03) - 2.03^2}{2.01} + \dfrac{s_U+k}{201}\\
&=\; \Omega(1)+\dfrac{s_U+k}{201},
\end{align*}
where the last equality holds since $k\ge 3$ and $s_U \ge 7$. Hence, $\eqref{eq union bd}\le \eqref{eq union bound 1}\le \left(\tfrac{k+s_U}{n}\right)^{\Omega(s_U+k)}$, so one may deduce by summing~\eqref{eq union bd} (or~\eqref{eq union bound 1}) over $s_U\in [s_{\min}, \eps n]$ that whp there is no set of $k-1$ vertices in $H_k$ separating a set of between 7 and $\eps n$ vertices from the rest.
    
The case when $s_U = \Omega(n)$ is also easily treated since the left hand side of~\eqref{eq union bd} is bounded from above by
\begin{align}
&\exp(o(n))\exp\left(s_U\log\left(\tfrac{n}{s_U}\right)+(n-s_U)\log\left(\tfrac{n}{n-s_U}\right) - \tfrac{s_Uk}{2}\log\left(\tfrac{n}{s_U}\right)\right)\nonumber\\ =\; 
&\exp(o(n))\exp\left((n-s_U)\log\left(\tfrac{n}{n-s_U}\right) - \tfrac{s_U(k-2)}{2}\log\left(\tfrac{n}{s_U}\right)\right),\label{eq case 1}
\end{align}
and since the function $t\in (0,1)\mapsto t\log(t^{-1})-(1-t)\log((1-t)^{-1})$ is non-negative on the interval $[0, 1/2]$ (note that, indeed, $s_U/n\in [0, 1/2]$),~\eqref{eq case 1} is $o(n^{-1})$ for $k\ge 5$. To extend the conclusion to the case when $k\in \{3,4\}$, a more precise upper bound than~\eqref{eq union bd} is needed. Such a bound can be provided by taking into account that, when one half-edge incident to a vertex in $U$ is matched, the number of available half-edges decreases by 2. Thus, one may replace~\eqref{eq union bd} by
\begin{align*}
\tbinom{n}{s_U,k-1} \left(\prod_{j=0}^{\lfloor s_U/2\rfloor}\tfrac{s_U+k-1-2j}{n-1-2j}\right)^{k}
&\le\hspace{0.3em}
\exp(o(n)) 
\left(\tfrac{n}{s_U}\right)^{s_U}
\left(\tfrac{n}{n-s_U}\right)^{n-s_U}
\left(\tfrac{s_U!! (n-s_U)!!}{n!!}\right)^{k}\\
&\le\hspace{0.3em}
\exp(o(n)) 
\left(\tfrac{n}{s_U}\right)^{s_U}
\left(\tfrac{n}{n-s_U}\right)^{n-s_U}
\left(\tfrac{s_U}{n}\right)^{s_Uk/2} \left(\tfrac{n-s_U}{n}\right)^{(n-s_U)k/2},
\end{align*}
which is $o(n^{-1})$ when $k\ge 3$.

We deal with the few remaining cases ``by hand''. Let us begin with the case $k\ge 3$ and $s_U = 2$. A direct first moment computation yields a bound of
\begin{equation*}
    \binom{n}{k+1} k^{2k} \frac{1}{(n-1)^k} \frac{1}{(n-3)^k}\le \frac{(ne)^{k+1}}{k^k} k^{2k} \frac{O(1)}{n^{2k}} = kO\left(\left(\frac{ek}{n}\right)^{k-1}\right).
\end{equation*}
Since we treat only the case $k = o(n)$, the above expression is always $o(n)$ (which can be seen, for example, by considering apart the cases $k\ge n^{1/3}$ and $k\le n^{1/3}$). We are left with a finite number of cases to check, namely $k\in [3,6]$ and $s_{U}\in [2, s_{\min}-1]$. Note that, in the case of a graph with $a = O(1)$ vertices and $b$ edges, the expected number of copies of this graph in $H_k$ is of order $\Theta(n^{a-b/2})$. Moreover, for a graph with at least $\lceil \tfrac{s_U k}{2}\rceil$ edges and $s_U+k-1$ vertices, one may easily check that $\lceil \tfrac{s_U k}{2}\rceil > s_U+k-1$ unless $(k, s_U)\in \{(3,3), (3,4), (4,3)\}$, hence the first moment computation ensures that whp there are no such subgraphs of $H_k$. To deal with the three remaining cases, note that, in each of them, $W$ must contain at least one vertex that is not connected to $U$ in $H_k$. Hence, ignoring isolated vertices in the first moment computation and using that $\lceil \tfrac{s_U k}{2}\rceil > s_U+k-2$ ensures that whp for all $k\ge 3$ and $s_U\in [2, n/2]$, a set of $s_U$ vertices in $H_k$ has at least $k$ neighbours.

Now, to make the graph $H_k$ $k$-connected whp, we must complete it to a random graph of minimum degree $k$. To do this, we must compensate the vertices which were incident to repeated edges. Denote by $Y$ the number of pairs of repeated edges (note that, if an edge appears in exactly $t$ different colors, it is counted $\tbinom{t}{2}$ times). Since $k = o(n)$, a union bound over all pairs of vertices in $V_n$ and all pairs of colours implies $\mathbb E [Y]\le \tbinom{n}{2} \tbinom{k}{2} \tfrac{1}{(n-1)^2}\le  \tfrac{k^2}{4}$, so by Markov's inequality $Y\le \frac{1}{2}\omega k^2$ whp. Conditioning on this event and using that $k = o(n)$, applying Lemma~\ref{lem bipart AB} for $m = \ell = n-k$ and a single vertex in the other part ensures that one can whp complete the graph $H_k$ to a graph of minimum degree $k$ by systematically adding an edge between a vertex of degree at most $k-1$ and a vertex of degree $k$ for $\omega k^2$ more rounds, which completes the proof of the first statement.\\

Now, if $k = o\left(\tfrac{n^{1/2}}{(\log n)^2}\right)$, we denote by $U_{\le r} = U_{\le r, i}$ the set of vertices of degree at most $r$ in $H_k$ at round $i$ (after identifying repeated edges in $H_k$) and $U_r = U_{r,i}$ the set of vertices of degree $r$ in $H_k$ at round $i$. 

We first treat the values of $k$ in the range $[\log\log n]$. Then, $\mathbb E [Y] = \Theta(k^2) = o\left((\log\log n)^2\right)$, so $Y = o\left(\log n\right)$ whp. Hence, it is sufficient to consecutively add additional edges to $H_k$ between $U_{\le k-1}$ to $U_k$ until $U_{\le k-1}$ becomes empty: indeed, at any step there is probability at least $1 - o\left(\tfrac{\log n}{n}\right)$ that one is able to construct such an edge. A union bound over at most $2Y = o(\log n)$ rounds shows that one obtains a graph with vertex degrees only $k$ and $k+1$ with probability $1-o\left(\tfrac{\log n}{n}\right) = 1-o(1)$.

Now, assuming that $k\ge \log\log n$, both $\mathbb E [Y] = \Theta(k^2)$ and, by an immediate application of the second moment method, $Y = (1+o(1))\mathbb E [Y]$ whp. Fix $\omega = \omega(n) = \log n$. Then, the number of vertices $Z_{\omega}$ incident to at least $\omega$ pairs of repeated edges satisfies 
$$\mathbb E [Z_{\omega}] \le (1+o(1)) n^{\omega+1} \tbinom{k}{2\omega} \tfrac{1}{(n-1)^{2\omega}} = O\left(n\left(\tfrac{k^2}{n-1}\right)^{\omega}\right) = o(1),$$
so $Z_{\omega}=0$ whp. Also, a similar computation shows that the number of vertices $Z_2$, incident to at least two pairs of repeated edges, satisfies $\mathbb E [Z_2] = O(\tfrac{k^4}{n}) = o(\tfrac{k^2}{\omega^3})$, so $Z_2 = o(\omega^{-3} Y)$ whp. We conclude that the pairs of repeated edges, incident to vertices, counted by $Z_2$, is at most $\omega^{-3} Y$ whp. Let us condition on all events in this paragraph that happen whp. We consider three regimes.

\paragraph{Regime 1.} As long as $U_{\le k-2}$ is not empty, look for an edge that connects a vertex in $U_{\le k-2}$ with a vertex in $U_k$. Since $|U_k| = n - o(n)$ throughout the process, every such attempt succeeds with probability $1-o(1)$. If the attempt at a given round is not successful, then try to connect a vertex in $U_{k-1}$ to a vertex in $U_k$: indeed, by Lemma~\ref{lem bipart AB} this is possible with probability $1 - \exp(-c |U_{k-1}|) = 1 - \exp(-\Omega(k^2))$. Thus, based on the fact that $\omega Z_2\le \omega^{-2} Y$, Markov's inequality guarantees that, within $\omega^{-1} Y = o(k^2)$ rounds, $U_{\le k-2}$ is emptied whp. Moreover, the application of Lemma~\ref{lem bipart AB} and a union bound ensure that no vertex of degree $k+1$ appears meanwhile whp: indeed, the probability for this event is $1 - O(k^2 \exp(-\Omega(k^2)) = 1-o(1)$.

\paragraph{Regime 2.} While $|U_{k-1}|\ge \omega$, an application of Lemma~\ref{lem bipart AB} for the complement of the bipartite graph $H[U_{k-1}, U_k]$ (note that $|U_k| = n-o(n)$ throughout the process) ensures that, at each step, there is probability at least $1 - \exp(-c |U_{k-1}|)$ to connect a vertex in $U_{k-1}$ to one in $U_k$. Hence, the probability of a failure at some step in this regime is at most $\sum_{j\ge \omega} \exp(-cj) = o(1)$.

\paragraph{Regime 3.} When $|U_{k-1}| < \omega$, we only use the fact that, at round $i$, any vertex $v\in V_n\setminus U_k$ connects to $U_k$ in $T_i$ with probability $1 - o(Y/n) = 1 - o(\omega^{-1})$ since $|U_k| = n-o(n)$ throughout the process. Thus, a union bound implies that $|U_{k-1}|$ goes from $\omega$ to 0 in exactly $\omega$ steps and no vertex of degree $k+2$ appears meanwhile whp, which completes the proof.
\end{proof}

\subsection{\texorpdfstring{Proof of Theorem~\ref{thm HC}}{}}

Our next goal is to show Theorem~\ref{thm HC}. Its proof will rest upon two main ingredients: Lemma~\ref{lem bipartite} and the following Lemma~\ref{lem long cycle}, which provides structural information for the union of two independent perfect matchings of $K_{n/2, n/2}$. In the next lemma, we assume that $n$ is even and denote $V_{\mathrm{even}} = \{v_{2i}\}_{i=1}^{n/2}$ and $V_{\mathrm{odd}} = V_n\setminus V_{\mathrm{even}}$. Despite the fact that the statement and its proof are classical for uniform random permutations, we provide a proof for completeness. Note that, since we work with simple graphs, a cycle of length 2 in our case is not a real cycle but only an edge that appears in both $M_1$ and $M_2$.

\begin{lemma}\label{lem long cycle}
Fix $\omega = \omega(n)\to \infty$. Let $M_1$ and $M_2$ be two independent perfect matchings between $V_{\mathrm{even}}$ and $V_{\mathrm{odd}}$. Then, whp $M_1\cup M_2$ contains at most $\omega \log n$ cycles and the largest of them is of length at least $\omega^{-1} n$.
\end{lemma}
\begin{proof}
For every even $i\in [n]$, denote by $W_i$ the number of cycles of length $i$ in $M_1\cup M_2$; then, 
\begin{equation*}
\mathbb E W_i = \binom{n/2}{i/2}^2 \left(\frac{i}{2}-1\right)!\left(\frac{i}{2}\right)! \prod_{j=0}^{i/2-1} \frac{1}{(n/2-j)^2} = \frac{2}{i},
\end{equation*}
where the binomial factors stand for the number of choices of $i$ vertices to constitute a cycle of length $i$, the term $\left(\tfrac{i}{2}-1\right)!\left(\tfrac{i}{2}\right)!$ stands for the number of different cycles one may obtain from $i$ vertices while respecting the parity constraint, and finally $\prod_{j=0}^{i/2-1} \frac{1}{(n/2-j)^2}$ stands for the probability that a fixed cycle of length $i$ indeed appears in $M_1\cup M_2$. Since
\begin{equation*}
\sum_{j=1}^{n/2} \mathbb E W_{2j} = \sum_{j=1}^{n/2} \frac{2}{2j} = (1+o(1))\log n \text{ and } \sum_{j=1}^{\omega^{-1} n/2} 2j\mathbb E W_{2j} = \omega^{-1}n = o(n),
\end{equation*}
we conclude by Markov's inequality that whp there are no more than $\omega \log n$ cycles, and moreover only at most $\omega^{-1/2} n$ vertices participate in cycles of length at most $\omega^{-1} n$, which proves the lemma.
\end{proof}

\begin{proof}[Proof of Theorem~\ref{thm HC}]
Fix $\hat\omega = \hat\omega(n) = \log n$ and $\omega = \omega(n)\to \infty$ satisfying $\omega(n) = n^{o(1)}$. We adopt the following strategy: first, we construct two independent perfect matchings $M_1$ and $M_2$ between $V_{\mathrm{even}}$ and $V_{\mathrm{odd}}$, and then we consecutively attempt to merge the cycles in $M_1\cup M_2$. By Lemma~\ref{lem bipartite} the first stage of the strategy takes at most $n + \omega^{1/2} n^{3/4}$ rounds whp, and by Lemma~\ref{lem long cycle} there is a cycle of length at least $\hat{\omega}^{-1} n$ and also at most $\hat{\omega} \log n$ cycles in total in $M_1\cup M_2$ whp. We condition on each of these events.

Now, we formally describe the second stage of the strategy. Let $(\pi_j)_{j=1}^k$ be all cycles in $M_1\cup M_2$, where $k\le \hat \omega \log n$ and $\pi_1$ is the longest cycle. Now, we define inductively a sequence of cycles $(\pi'_j)_{j=1}^k$ as follows. First, set $\pi'_1 = \pi_1$ and suppose that, for some $j\in [2,k]$, $\pi'_{j-1}$ has been defined at round $i_{j-1}$ (in particular, $i_1-1$ is the last round of the construction of $M_1\cup M_2$). Then, fix a vertex $u_j\in \pi_j$, a neighbour $w_j$ of $u_j$ in $\pi_j$, and denote by $N_j(t)$ the set of vertices in $\pi'_{j-1}$, sharing a common neighbour in $\pi'_{j-1}$ with $u_j$ after round $i_{j-1}+t$, see Figure~\ref{fig 3}. 

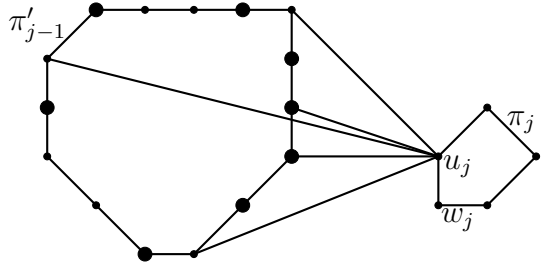
\begin{figure}
\centering
\begin{tikzpicture}[scale = 1.3,line cap=round,line join=round,x=1cm,y=1cm]
\clip(-5.771256710309041,-0.6) rectangle (5,2.2);
\draw [line width=0.8pt] (-2,2)-- (-2.5,2);
\draw [line width=0.8pt] (-2.5,2)-- (-3,1.5);
\draw [line width=0.8pt] (-3,1.5)-- (-3,1);
\draw [line width=0.8pt] (-3,1)-- (-3,0.5);
\draw [line width=0.8pt] (-3,0.5)-- (-2.5,0);
\draw [line width=0.8pt] (-2.5,0)-- (-2,-0.5);
\draw [line width=0.8pt] (-2,-0.5)-- (-1.5,-0.5);
\draw [line width=0.8pt] (-1.5,-0.5)-- (-1,0);
\draw [line width=0.8pt] (-1,0)-- (-0.5,0.5);
\draw [line width=0.8pt] (-0.5,0.5)-- (-0.5,1);
\draw [line width=0.8pt] (-0.5,1)-- (-0.5,1.5);
\draw [line width=0.8pt] (-0.5,1.5)-- (-0.5,2);
\draw [line width=0.8pt] (-0.5,2)-- (-1,2);
\draw [line width=0.8pt] (-1,2)-- (-1.5,2);
\draw [line width=0.8pt] (-1.5,2)-- (-2,2);
\draw [line width=0.8pt] (1,0.5)-- (1.5,1);
\draw [line width=0.8pt] (1.5,1)-- (2,0.5);
\draw [line width=0.8pt] (2,0.5)-- (1.5,0);
\draw [line width=0.8pt] (1.5,0)-- (1,0);
\draw [line width=0.8pt] (1,0)-- (1,0.5);
\draw [line width=0.8pt] (1,0.5)-- (-0.5,0.5);
\draw [line width=0.8pt] (1,0.5)-- (-0.5,2);
\draw [line width=0.8pt] (1,0.5)-- (-1.5,-0.5);
\draw [line width=0.8pt] (1,0.5)-- (-3,1.5);
\draw [line width=0.8pt] (1,0.5)-- (-0.5,1);
\begin{scriptsize}
\draw [fill=black] (-2,2) circle (1pt);
\draw [fill=black] (-2.5,2) circle (2pt);
\draw [fill=black] (-3,1.5) circle (1pt);
\draw [fill=black] (-3,1) circle (2pt);
\draw [fill=black] (-3,0.5) circle (1pt);
\draw [fill=black] (-2.5,0) circle (1pt);
\draw [fill=black] (-2,-0.5) circle (2pt);
\draw [fill=black] (-1.5,-0.5) circle (1pt);
\draw [fill=black] (-1,0) circle (2pt);
\draw [fill=black] (-0.5,0.5) circle (2pt);
\draw [fill=black] (-0.5,1) circle (2pt);
\draw [fill=black] (-0.5,1.5) circle (2pt);
\draw [fill=black] (-0.5,2) circle (1pt);
\draw [fill=black] (-1,2) circle (2pt);
\draw [fill=black] (-1.5,2) circle (1pt);
\draw [fill=black] (1,0.5) circle (1pt);
\draw[color=black] (1.2,0.4) node {\large{$u_j$}};
\draw[color=black] (1.85,0.85) node {\large{$\pi_j$}};
\draw[color=black] (-3.1,1.85) node {\large{$\pi'_{j-1}$}};
\draw [fill=black] (1.5,1) circle (1pt);
\draw [fill=black] (2,0.5) circle (1pt);
\draw [fill=black] (1.5,0) circle (1pt);
\draw [fill=black] (1,0) circle (1pt);
\draw[color=black] (1.2,-0.15) node {\large{$w_j$}};
\end{scriptsize}
\end{tikzpicture}
\caption{The graph after round $i_{j-1}+t$; the vertices in $N_j(t)$ are thickened. Note that, for every $w'\in N_j(t)$, connecting $w_j$ to $w'$ ensures the existence of a cycle, containing both $\pi'_{j-1}$ and $\pi_j$; this is exactly the cycle $\pi'_j$.}
\label{fig 3}
\end{figure}

Then, at any round $i_{j-1} + t$ until $\pi'_j$ is not defined, do the following:
\begin{itemize}
    \item if the tree, proposed at round $i_{j-1}+t$, connects $w_j$ to $w'_j\in N_j(t)$, construct the edge $w_jw'_j$. Then, find a vertex $u'_j\in \pi'_{j-1}$, connected to both $w'_j$ and $u_j$, and define the cycle $\pi'_j$ as the union of the edges $u_ju'_j$, $w_jw'_j$, $\pi'_{i-1}\setminus \{u'_jw'_j\}$ and $\pi_j\setminus \{u_jw_j\}$. Also, define $i_j = i_{j-1}+t$;
    \item if the previous point does not hold and the tree, proposed at round $i_{j-1}+t$, connects $u_j$ to a non-neighbour (in $G_{i_{j-1}+t-1}$) $u'\in \pi'_{j-1}$, set $G_{i_{j-1}+t} = G_{i_{j-1}+t-1}\cup \{u_ju'\}$;
    \item else, skip the round.
\end{itemize}
Note that, since there are always at least $\hat{\omega}^{-1} n$ vertices in $\pi'_{j-1}$ ($\pi'_{j-1}$ has at least as many vertices as $\pi_1$), for any of the first $\hat{\omega}^{-1} n/2$ rounds after $i_{j-1}$, the second point of the algorithm ensures that there is probability at least $(2\hat \omega)^{-1}$ that the round is not skipped.

For all $j\in [2, k]$, denote by $\cE_j(t)$ the event $\{i_j - i_{j-1}\le t\}$, and also denote by $\cE_j'(t)$ the event ``$u_j$ has at least $\hat \omega n^{1/2}$ neighbours in $\pi'_{j-1}$ after round $i_{j-1}+t$''. On the one hand, Chernoff's inequality and the fact that $\hat{\omega}^3 n^{1/2}\ll \hat{\omega}^{-1} n/2$ ensure that $\Prob(\cE_j'(\hat{\omega}^3 n^{1/2})\mid \overline{\cE_j(\hat{\omega}^3 n^{1/2})}) = 1-o(n^{-1})$. On the other hand, conditionally on $\cE_j'(\hat{\omega}^3 n^{1/2})$, $N_j(i_{j-1}+\hat{\omega}^3 n^{1/2})\ge \hat{\omega} n^{1/2}$. Now, for every $t\ge 1$, denote by $X_t$ the indicator random variable of the event ``$w_i$ connects to a vertex in $N_j(\hat{\omega}^3 n^{1/2})$ at round $i_{j-1}+\hat{\omega}^3 n^{1/2}+t$''. Then, conditionally on $\cE_j'(\hat{\omega}^3 n^{1/2})$ (and also on $\cE_j'(\hat{\omega}^3 n^{1/2})\cap \overline{\cE_j(\hat{\omega}^3 n^{1/2})}$), $\Prob(X_t = 1)\ge \hat{\omega} n^{-1/2}$ for all $t\ge 1$, so by Chernoff's inequality
\begin{align*}
\Prob(\cE_j(2\hat{\omega}^3 n^{1/2}))
&\ge \Prob(\cE_j(2\hat{\omega}^3 n^{1/2})\mid \overline{\cE_j(\hat{\omega}^3 n^{1/2})})\\
&\ge \Prob(\cE_j(2\hat{\omega}^3 n^{1/2})\cap \cE_j'(\hat{\omega}^3 n^{1/2})\mid \overline{\cE_j(\hat{\omega}^3 n^{1/2})})\\
&\ge \Prob(\cE_j(2\hat{\omega}^3 n^{1/2})\mid \cE_j'(\hat{\omega}^3 n^{1/2})\cap \overline{\cE_j(\hat{\omega}^3 n^{1/2})}) - \Prob(\overline{\cE_j'(\hat{\omega}^3 n^{1/2})}\mid \overline{\cE_j(\hat{\omega}^3 n^{1/2})})\\
&\ge \Prob\left(\sum_{t=1}^{\hat{\omega}^3 n^{1/2}} X_t = 0\;\bigg|\; \cE_j'(\hat{\omega}^3 n^{1/2})\cap \overline{\cE_j(\hat{\omega}^3 n^{1/2})}\right) - o(n^{-1})\\ 
&= 1-o(\exp(-\Omega(\hat{\omega}^4))+n^{-1}) = 1-o(n^{-1}).
\end{align*}
A union bound over all $k-1$ cycle mergings shows that whp the second stage takes at most $2\hat{\omega}^4 n^{1/2} \log n$ rounds, which leads to at most $n + \omega^{1/2}n^{3/4} + 2\hat{\omega}^4 n^{1/2} \log n\le n + \omega n^{3/4}$ rounds in total whp.
\end{proof}

\subsection{\texorpdfstring{Proof of Theorem~\ref{thm factors}}{}}

The following lines prove Theorem~\ref{thm factors}. We remark that, roughly speaking,~\eqref{pt 1f} and~\eqref{pt 2f} in Theorem~\ref{thm factors} rely on the fact that we are able to construct perfect matchings fast thanks to Lemma~\ref{lem bipartite}.

\begin{proof}[Proof of Theorem~\ref{thm factors}~\eqref{pt 1f}]
Set $k = |H|$ and fix an arbitrary vertex $r$ of $H$ to serve as a root. Then, orient all edges away from $r$ and perform a depth-first search (DFS) of $H$. Let $(e_1,e_2,\dots,e_{k-1})$ be the list of edges of $H$ in order of their appearance in the DFS. For all $i\in [k-1]$, denote by $e^-_i$ the vertex in $e_i$ closer to $r$, and by $e^+_i$ the vertex in $e_i$ further from $r$. 

Now, we describe the strategy $\sigma$ that we will use to form an $H$-factor. It is divided into $k-1$ stages, which are similar in nature. In the beginning, fix a set $U(e^-_1)\subseteq V_n$ of size $n/k$. For every $i\in [k-1]$, at the $i$-th stage we construct a matching between $U(e^-_i)$ and $V_n\setminus \cup_{j=1}^{i} U(e^-_j)$ of size $n/k$, and define the set $U(e^+_i)$ as the part of the matching, contained in $V_n\setminus \cup_{j=1}^{i} U(e^-_j)$. By Lemma~\ref{lem bipartite} this construction takes at most $n/k + O(\omega^{1/2} n^{3/4})$ rounds for any $\omega = \omega(n)\to \infty$ whp. Hence, by a union bound, the number of rounds needed to complete all $k-1$ stages is at most $\tfrac{(k-1)n}{k} + \omega n^{3/4}$ whp.
\end{proof}

We turn our attention to the proof of Theorem~\ref{thm factors}~\eqref{pt 2f}. 

\begin{lemma}\label{lem almost factor}
Given any graph $G$ (possibly disconnected) and any $\eps > 0$, there is a constant $C = C(G, \varepsilon) > 0$ and a strategy $\sigma\in \mathcal S_n$ that constructs $(1-\varepsilon) \tfrac{n}{|G|}$ disjoint copies of $G$ in at most $C n$ rounds whp.
\end{lemma}
\begin{proof}
Set $\delta = \tfrac{\varepsilon}{|E(G)|}$, $s = |G|$, and partition $V_n$ into $\tfrac{n}{s}$ sets $(W_i)_{i=1}^{n/s}$ of $s$ vertices. Then, for any set $W_i$, choose an arbitrary bijection $\phi_i: v\in V(G)\mapsto w\in W_i$. 

We divide the construction into $|E(G)|$ stages as follows. At any stage, fix an edge $uv\in E(G)$ and proceed by attempting to construct the edges $E_{uv} = \{\phi_i(uv)\}_{i=1}^{n/s}$. If one cannot construct a new edge at a given round, then this round is ignored. Let us show that there is a constant $C_1 = C_1(G, \delta) > 0$ satisfying that whp one needs at most $C_1 n$ steps to construct a $(1-\delta)$-proportion of $E_{uv}$: indeed, this is sufficient for our purposes since after the end of the last stage there are at least $(1-|E(G)|\delta) \tfrac{n}{s} = (1-\eps) \tfrac{n}{s}$ copies of $G$ in the USTSR graph.

For every $j\in [(1-\delta) n/s]$, let $Z_j$ be the time needed to construct the $j$-th edge in $E_{uv}$. Then, by Lemma~\ref{lem independence} we deduce that $Z_j$ is a geometric random variable with parameter $1 - \left(1-\tfrac{2}{n}\right)^{n/s - i} = 1 - \exp(-(1+o(1)) \tfrac{2}{s} - \tfrac{2i}{n})$, which is minimal when $i = (1-\delta) \tfrac{n}{s}$, and in this case this gives $1 - \exp(-(1+o(1)) \tfrac{2\delta}{s})\ge p := 1 - \exp(-\tfrac{\delta}{s})$. Then, for every sufficiently small $\lambda > 0$ one has
\begin{align*}
\Prob\left(\sum_{j=1}^{(1-\delta) n/s} Z_j\ge \frac{2n}{p} \right)
&\le\hspace{0.3em} 
\frac{\mathbb E[\exp(\sum_{j=1}^{(1-\delta) n/s} \lambda Z_j)]}{\exp(2\lambda n/p)}\\ 
&=\hspace{0.3em} 
\frac{\prod_{j=1}^{(1-\delta) n/s} \mathbb E[\exp(\lambda Z_j)]}{\exp(2\lambda n/p)}\\
&=\hspace{0.3em} 
\frac{\left(\mathbb E[\exp(\lambda \mathrm{Geom}(p))]\right)^{n}}{\exp(2\lambda n/p)}\\
&\le\hspace{0.3em} \left(\frac{(1-p)\exp(\lambda)}{p(1-p\exp(\lambda))\exp(2\lambda/p)}\right)^n = o(1).
\end{align*}
Hence, each of the $|E(G)|$ stages requires at most $\tfrac{2}{p} n$ rounds to be completed whp, which proves the lemma for $C_1 = \tfrac{2}{p}$.
\end{proof}

\begin{proof}[Proof of Theorem~\ref{thm factors}~\eqref{pt 2f}]
Recall that, for a central graph $H$, we call \emph{witness} of $H$ any vertex $u\in V(H)$ such that every edge $e\in E(H)$, incident to $u$, satisfies that $H\setminus e$ is a disconnected graph. Fix $k = |H|$ and a witness $u$ of $H$ with neighbourhood $\{u_1,u_2,\dots,u_t\}$. Set $\sigma\in \mathcal S_n$ to be the strategy that, first, constructs $n/k$ copies of $G := H\setminus v$ following the strategy from Lemma~\ref{lem almost factor} (applied for $G$ and $\varepsilon = 1/k$), and then aims to add $t = \deg_H(u)$ perfect matchings to extend the already present disjoint copies of $G$ to an $H$-factor. More formally, for all $i\in [r]$, let $U_i$ be the set of $\tfrac{n}{k}$ vertices, corresponding to the vertex $u_i$ in the $\tfrac{n}{k}$ copies of $G$, and let $U$ be the set of vertices not belonging to any copy of $G$. After generating the $\tfrac{n}{k}$ copies of $G$ in $C(G, \eps) n$ steps whp, it remains to construct $t$ matchings of size $\tfrac{n}{k}$ between $U_0$ and $U_i$ for all $i\in [t]$. By Lemma~\ref{lem bipartite} applied $t$ times we deduce that any $C_H > C(G, \eps) + \tfrac{t}{k}$ satisfies the requirement.
\end{proof}

Before proving Theorem~\ref{thm factors}~\eqref{pt 3f}, let us provide the rough idea behind our argument. Fix a set of $r$ vertices labeled $1,2, \ldots, r$ and set $R = \tbinom{r}{2}$. For every $j\in [r-1]$, let $H_j$ be a path on $j+1$ vertices starting at 1 and ending at $r$, and let $H_{r-1}\subseteq H_r\subseteq \dots \subseteq H_R = K_r$ be a sequence of graphs on the same vertex set such that, for every $i\in [R]$, $|E(H_i)| = i$. Using the differential equations method (Theorem~\ref{Thm:DEMethod}), we will try to build multiple copies of each of $H_1, \ldots, H_R$ and then connect them to form $\Omega(n)$ long chains. The number of copies of $H_i$ in every chain will be chosen the same and will depend (in some implicit way) on the number of rounds needed to complete the necessary number of copies of $H_{i+1}, \dots, H_R$. However, this number of rounds happens to be well concentrated by Theorem~\ref{Thm:DEMethod}.

We now concentrate on the formal description of the proof.

\begin{proof}[Proof of Theorem~\ref{thm factors}~\eqref{pt 3f}]
Our strategy of choice is divided into two stages. At the first stage, we recursively construct disjoint subsets $U_1', U_2', \ldots$ of $V_n$ of size $r$ as follows. In the beginning, one has $\nu(0) = 0$ sets. Suppose that, right before round $t\ge 1$, we have constructed $\nu(t-1)$ sets $U_1', U_2', \ldots, U_{\nu(t-1)}$. Also, fix positive integers $(\ell_i)_{i=0}^{\binom{r}{2}}$ to be chosen in a bit, and do the following: at round $t$, let $i\in [0, R]$ be the largest integer satisfying that the number of sets among $(U_j')_{j=1}^{\nu(t-1)}$ satisfying $G[U_j'] = H_i$ is less than $\tfrac{n}{\ell_i r}$. Then, set $i'\leftarrow i-1$ and, while $i'\neq 0$, look for an edge in $T_t$ that completes a copy of $H_{i'}$ among $(G[U_j'])_{j=1}^{\nu(t-1)}$ to form a copy of $H_{i'+1}$. If the attempt is successful, add an edge to the copy $G[U'_j]$ of $H_{i'}$ with the smallest possible $j$ and proceed to the next round. If the attempt is not successful, do $i'\leftarrow i'-1$ and come back to the loop. Then, either the round is completed by adding an edge within some of the sets $(U_j')_{j=1}^{\nu(t-1)}$, or we end up with $i'=0$. In the second case, connect two vertices $u,v\in V_n\setminus (U_j')_{j=1}^{\nu(t-1)}$, then pick $r-2$ arbitrary vertices $w_1,\ldots,w_{r-2}$ in $V_n\setminus ((U_j')_{j=1}^{\nu(t-1)}\cup \{u,v\})$, and form a new set $U'_{\nu(t-1)+1} = \{u,v,w_1,\ldots,w_{r-2}\}$ inducing a single edge (so a new copy of $H_1$, up to giving labels 1 and $r$ to $u$ and $v$). If even this is not possible, stop the process and end the first stage.

Denote by $X_i(t)$ the number of sets among $(U_j')_{j=1}^{\nu(t)}$ that satisfy $G[U_j'] = H_i$ after round $t$, and denote by convention $X_{R+1}\equiv 0$. We will show via the differential equation method that one may choose the integers $\ell_1, \dots, \ell_R$ so that the first stage ends at round $t_f\le n/4$ whp and moreover, for every $i\in [R]$, $X_i(t_f) = \tfrac{n}{\ell_i r}$. First, note that, at any round $t\ge 0$, each of the numbers $(|X_i(t) - X_i(t+1)|)_{i=1}^{R}$ is 0 or 1. Moreover, for any $t\ge 0$, conditionally on $(X_i(t))_{i=1}^{R}$, for any $j\in [2,R]$ the expected value of $X_j(t+1)$ is given by 
$$X_j(t) - \left(1-\tfrac{2}{n}\right)^{X_j(t)}\mathds{1}_{j < R} + \left(1 - \left(1-\tfrac{2}{n}\right)^{X_j(t)}\mathds{1}_{j < R}\right)\left(1-\tfrac{2}{n}\right)^{X_{j-1}(t)},$$ 
and a similar expression may be derived for $X_1(t)$. Indeed, to form a new copy of $H_j$ one has to hit one of $X_{j-1}(t)$ disjoint non-edges in the copies of $H_{j-1}$, and, except if $j = R$, in order to transform one copy of $H_j$ into $H_{j+1}$ one has to hit one of $X_j(t)$ disjoint non-edges in the copies of $H_j$ (by definition, priority is given to the second operation until the necessary number of copies of $H_{j+1}$ is attained). Therefore, the differential equation method (Theorem~\ref{Thm:DEMethod}) guarantees that the trajectories of $(X_i(t))_{i=0}^{R}$ are concentrated around their expected values with exponentially small deviation probability. Hence, one may sequentially choose $\ell_1, \ldots, \ell_{R}$ so that, for every $i\in [R]$, $\ell_i$ is so large with respect to $(\ell_j)_{j=1}^{i-1}$ that, when $X_i(t) = \tfrac{n}{\ell_i r}$ for the first time at some round $t$, whp $X_j(t) < \tfrac{n}{\ell_i r}$ for all $j\in [i-1]$. Apart from that, choosing $\ell_1$ sufficiently large ensures that, at round $n/4$, say, $X_i(t) = \tfrac{n}{\ell_i r}$ for all $i\in [R]$ whp. Thus, we choose $\ell_1, \ldots, \ell_{\tbinom{r}{2}}$ so that the above requirements are met and condition on the success of this first stage.

Once the required number of copies of $H_1, \ldots, H_{\tbinom{r}{2}}$ have been constructed, divide the set of connected components into equal sets $S_1,\ldots,S_m$, where $m\in \mathbb N$ is fixed, and for every $i\in [m]$, $S_i$ contains copies of only one graph $H'_i$ among $H_0$, which corresponds to an isolated vertex, and $H_1, \ldots, H_{\tbinom{r}{2}}$. Then, $m-1$ consecutive applications of Lemma~\ref{lem bipartite}, each providing a perfect matching between $\Omega(n)$ copies of the same vertex in $H'_i$ within $S_i$ and $\Omega(n)$ copies of the same vertex in $H'_i$ within $S_{i+1}$, construct whp copies of the same graph $H$ (which, roughly speaking, resembles a chain connecting many copies of $H_1, \ldots, H_R$) within $\tfrac{n |E(H)|}{|H|} + m \omega n^{3/4}$ rounds for any function $\omega = \omega(n)\to \infty$. Since the graph $H$ contains at least one copy of $H_R = K_r$, the strategy $\sigma$ satisfies the statement of the theorem for $H$.
\end{proof}

\subsection{\texorpdfstring{Proof of Theorem~\ref{thm covering}}{}}

Fix a graph $H$ and a vertex $u$ of $H$ of minimum degree, do a depth first search (DFS) from $u$ and order the vertices accordingly. Then, set $V(H) = \{u_0 = u, u_1, \ldots, u_k\}$ where $N_H(u) = \{u_1,\ldots,u_r\}$. Fix $\omega = \omega(n) = (\log n)^3$. The following algorithm constitutes an important part of the strategy $\sigma$ that will be analyzed in the proof of Theorem~\ref{thm covering}.

\begin{algorithm}\label{algo 2}
\emph{Input: an empty graph on $n$ vertices.
\noindent
Set $i\leftarrow k$. For all $j\in [k-1]$, set $U_j=\emptyset$ and $U_k \leftarrow \{v_{n-\omega+1},\ldots,v_n\}$. Set $E = \emptyset$. While $i\neq 0$, do:
\begin{enumerate}[(i)]
    \item\label{pt 2.1} Set $F = \{u_iu_j: 0<j<i \text{ and } u_iu_j\in E(H)\}$.
    \item\label{pt 2.2} If $|U_i| > \omega$, delete arbitrary $|U_i|-\omega$ elements from $U_i$.
    \item\label{pt 2.3} For every $j$ such that $u_iu_j\in F$ and for every $v\in U_i$, do:
        \begin{itemize}
            \item if $U_j=\emptyset$, for $2\omega^{-2} n$ rounds, attempt to construct an edge between $v$ and $V_n\setminus \cup_{t=1}^k U_t$. At every successful round, add the new neighbour of $v$ to $U_j$. Interrupt the algorithm if the vertex $v$ has had less than $\omega^{-2} n$ neighbours after $2\omega^{-2} n$ rounds.
            \item else, if $U_j\neq \emptyset$, for $2\omega^{-2} n$ rounds, attempt to construct a new edge between $v$ and $U_j$. Interrupt the algorithm if the vertex $v$ has had less than $\omega^{-2} |U_j|$ neighbours after $2\omega^{-2} n$ rounds. Finally, do $U_j\leftarrow U_j\cap N_{U_j}(v)$.
        \end{itemize}
    \item\label{pt 2.4} Do $i\leftarrow i-1$.
\end{enumerate}}
\end{algorithm}

\begin{claim}\label{cl algo 2}
Algorithm~\ref{algo 2} always terminates in $O(\omega^{-1} n)$ steps. Moreover, in the end of the algorithm, $i=0$ whp.
\end{claim}
\begin{proof}
The first statement is due to the fact that step~\eqref{pt 2.3} of Algorithm~\ref{algo 2} is iterated $\omega E(H)$ many times, and every iteration requires at most $2\omega^{-2} n$ rounds. For the second step, we show that step~\eqref{pt 2.3} of Algorithm~\ref{algo 2} interrupts the algorithm with probability only $o(\omega^{-1})$ for every $u_iu_j\in F$ and $v\in U_i$. We consider two cases:
\begin{itemize}
    \item if $U_j = \emptyset$, then at each of the following $2\omega^{-2}n$ rounds the vertex $v$ has probability $1-o(1)$ to be able to connect to a new neighbour since $|V_n\setminus \cup_{t=1}^k U_t| = n-o(n)$ by definition. By Chernoff's inequality, applied for the indicator random variables that $v$ connects to a new neighbour at one of the $2\omega^{-2}n$ given rounds, $\Prob(|U_j|\le \omega^{-2} n \text{ in the end})\le \exp(-\Omega(\omega^{-2} n)) = o(\omega^{-1})$;
    \item if $U_j \neq \emptyset$, then by construction $|U_j| = \omega^{-O(1)} n$. Then, while $v$ has less than $|U_j|/3$ neighbours, at each of the following $2\omega^{-2}n$ rounds the vertex $v$ has probability at least $\tfrac{2|U_j|}{3n}$ to be able to connect to a new neighbour in $U_j$. By Chernoff's inequality, applied to the indicator random variables that $v$ connects to a new neighbour in $U_j$ at one of the $2\omega^{-2}n$ given rounds, $\Prob(|N_{U_j}(v)\cap U_j|\le \omega^{-2} |U_j| \text{ in the end})\le \exp(-\Omega(\omega^{-2} |U_j|)) = o(\omega^{-1})$;
\end{itemize}
Thus, a union bound over the $E(H)$ edges of $H$ and the $\omega$ copies of every vertex $u_i\in H\setminus \{u\}$, left in $U_i$ after step~\eqref{pt 2.2} of the algorithm, proves the claim.
\end{proof}

After the end of Algorithm~\ref{algo 2}, fix $V'_n = V_n\setminus \cup_{i=1}^k U_i$. Also, define by $G$ the graph, constructed by Algorithm~\ref{algo 2} after it terminates.

\begin{claim}\label{cl 2 algo 2}
If $i=0$ in the end of Algorithm~\ref{algo 2}, then any graph $G'$, obtained by connecting each vertex in $V'_n$ with one vertex in each of $(U_i)_{i=1}^r$, contains a covering of $V_n$ by copies of $H$.
\end{claim}
\begin{proof}
For every $t\in [k]$, choose any vertex $w_t\in U_t$. Then, by construction of Algorithm~\ref{algo 2}, $w_t$ sends edges towards each of $\{w_j: 0<j<t \text{ and } u_tu_j\}$. Hence, if $i=0$ when Algorithm~\ref{algo 2} terminates, then the map $\varphi: u_t\in V(H)\setminus \{u_0\}\mapsto w_t$ is such that $\varphi(H\setminus \{u_0\})\subseteq G[\{w_t\}_{t=1}^k]$. Hence, connecting a new vertex to $\{w_t\}_{t=1}^r$ completes $G[\{w_t\}_{t=1}^k]$ to a copy of $H$.
\end{proof}

\begin{proof}[Proof of Theorem~\ref{thm covering}]
By Claim~\ref{cl algo 2} and Claim~\ref{cl 2 algo 2} it remains to show that one can match every vertex in $V'_n$ with each of $(U_t)_{t=1}^r$ in $d_{\min}(H) n + o(n)$ rounds. In fact, we show that one can match every vertex in $V'_n$ with a vertex in $U_1$ in $n+o(n)$ rounds whp; matching $V'_n$ with $(U_t)_{t=2}^k$ is done in the same way. The proof of this fact resembles the end of the proof of Lemma~\ref{lem bipartite}; indeed, for every $t\in [|V'_n|]$, define $Z_t$ to be the time one needs to reduce the number of unmatched vertices in $V'_n$ from $t$ to $t-1$. Then, $Z_t$ is a geometric random variable with parameter $\Prob(E(T)\cap (V_t\times (V_{t+\omega}\setminus V_t)) \neq \emptyset)$, where $T\sim \mathrm{Unif}(\cT_n)$, which is given by Lemma~\ref{lem 2.8}:
$$1 - \Prob(E(T)\cap (V_t\times (V_{t+\omega}\setminus V_t)) = \emptyset)\ge 1 - \left(1-\tfrac{t}{n}\right)^{\omega-1}\left(1-\tfrac{\omega}{n}\right)^{t-1}\ge 1-\exp\left(-\tfrac{2t\omega - t - \omega}{n}\right)\ge 1-\exp\left(-\tfrac{t\omega}{2n}\right).$$

\noindent
Hence, we get that
\begin{align*}
\sum_{t=1}^{|V'_n|} \mathbb E [Z_t-1] 
\le\; 
&\sum_{t=1}^{|V'_n|} \frac{1}{1-\exp\left(-\tfrac{t\omega}{2n}\right)} - 1\\
=\; 
&\sum_{t=1}^{|V'_n|} \frac{\exp\left(-\tfrac{t\omega}{2n}\right)}{1-\exp\left(-\tfrac{t\omega}{2n}\right)}\\ 
=\; 
&\sum_{t=1}^{\omega^{-1}n} \frac{\exp\left(-\tfrac{t\omega}{2n}\right)}{1-\exp\left(-\tfrac{t\omega}{2n}\right)} + \sum_{t=\omega^{-1}n+1}^{|V'_n|} \frac{\exp\left(-\tfrac{t\omega}{2n}\right)}{1-\exp\left(-\tfrac{t\omega}{2n}\right)}\\ 
=\; 
&\sum_{t=1}^{\omega^{-1}n} O(\tfrac{n}{t\omega}) + \sum_{t=\omega^{-1}n+1}^{|V'_n|} O(\exp\left(-\tfrac{t\omega}{2n}\right))\\
=\; 
&O(\tfrac{n\log n}{\omega}) = o(\tfrac{n}{\omega^{1/2}}).
\end{align*}

Thus, by Markov's inequality for the random variable $\sum_{t=1}^{|V'_n|} (Z_t-1)$, matching every vertex in $V'_n$ to a vertex in $U_1$ requires no more than $n + \omega^{-1/3} n$ steps whp. Repeating the same reasoning for each of $U_2, \ldots, U_r$ concludes the proof of the theorem.
\end{proof}

\bibliography{bib}
\bibliographystyle{plain}
\end{document}